\documentclass [12pt]{article}
\usepackage{graphicx,amssymb,amsfonts,latexsym,amsmath,amsthm,times}
\usepackage{epsfig}
\usepackage{fancyhdr}
\usepackage{color}
\usepackage[english]{babel}
\numberwithin{equation}{section}

\newtheorem{theorem}{Theorem}
\newtheorem{proposition}[theorem]{Proposition}
\newtheorem{remark}[theorem]{Remark}

\begin{document}

\author{D. Melchionda$^a$, E. Pastacaldi$^a$, C. Perri$^a$ and E. Venturino$^a$ \footnote{Corresponding
author. E-mail:
%diego.melchionda@gmail.com, elisa.pastacaldi@libero.it, chiarapperri@gmail.com, 
ezio.venturino@unito.it}\\
$^\dagger$Dipartimento di Matematica
via Carlo Alberto 10,\\
Universit\`a di Torino,\\
10123 Torino, Italy
}
\title{Interacting population models with pack behavior}
\date{}
\maketitle

\begin{abstract}
Models of coordinated behavior of populations living in the same environment are introduced for the cases when
they either compete with each other, or they both gain by mutual interactions, or finally when one hunts the other one.
The equilibria of the systems are analysed, showing that in some cases the populations may both disappear.
Coexistence leads to global asymptotic stability for symbiotic populations, or to Hopf bifurcations for
predator-prey systems. Finally, a new very interesting phenomenon is discovered in one of these rather simple
models. Indeed
tristability may be achieved in the competition case, for which competitive
exclusion is allowed to occur together with populations coexistence.
\end{abstract}



\section{Introduction}

We consider models for populations whose intermingling may be beneficial to both of them, harmful for both of them, or
beneficial for one and detrimental for the other one. The classical models always assume individualistic behavior of each
population, see e.g. Part I of \cite{MPV}. New models for mimicking the herd group defense of herbivores have been recently
introduced,
\cite{1}, differing quite a bit from older ideas relying on different assumptions, \cite{FW}, or from more
recent contributions, \cite{GG}, in which the word ``group'' is used with a completely different meaning,
or at least it is mathematically modeled in a completely different way from ours.
Furthermore, the biological literature abounds on social, herd or pack behaviour, but these concepts are either not the same that
are considered here or modeled via different mathematical tools, e.g. graph theory or game theory, than those employed here,
see for instance \cite{Giardina,Sumpter} and the wealth of literature that is cited in these papers.

Although \cite{1} deals with a demographic model, an extension to ecoepidemic systems has been proposed,
\cite{EV_JBS}. These kind of systems consider basic population interactions, on which the effects of
epidemics are superimposed. In the past twenty years a great deal of research effort has been devoted
to their understanding, for a brief introduction, see Chapter 7 of \cite{MPV}. These investigations represent
a natural outgrowth of the developments of the late nineties in mathematical epidemiology, when the effects of
population changes began to be accounted for in models for the spreading of diseases, \cite{B-VdD,G-H,ML-H,H00}.

In this paper we confine ourselves to the pure demographic situation, however. We
extend the community gathering idea also to predators, allowing them to hunt the prey in a coordinate fashion.
Thus the models we consider markedly differ from those of
\cite{1}, in the biological assumptions and above all in the mathematical structure.
In fact, this is not a mere extension of previous work to different population associations, because the predators are behaving individualistically
in
\cite{1} and therefore predator's pack behavior is absent in those models.

When the predators' pack
hunts the prey in general the individuals that have the major benefit are those that either take the most advantageous positions
in the community in order to get the best share of the loot, or simply those that get it because they are stronger.
Assume therefore that positions on the edge of the pack have the best returns for
the individuals that occupy them, since they are the first to fall upon the prey.
The main idea of community behaviors for predators had been considered in \cite{cosner}.
In the framework of animals socialized behavior
these ideas have recently been discussed also in \cite{4} and carried over to ecoepidemic systems, \cite{CV}.
We consider two situations for the prey, namely when they behave individualistically or when they gather in herds,
following the assumptions of
\cite{1}. In the latter situation, the most harmed prey during predators' hunting are those staying
on the boundary of the herd.

Here we also extend this concept to more general types of interactions among populations thriving in the same environment.
The cases of symbiosis and competition are also well-known in the literature. Again, the classical
approach envisions an individualistic behavior for the involved populations. In part this idea has been introduced in \cite{1},
but assuming that only one population behaves socially, the individuals of the other one live independently of each other.
We extend now the analysis to the case in which both populations show a community behavior, both when each one of the two communities benefits
from the interactions with the other one, as well as the case in which the communities compete with each other.

The systems introduced here mathematically model the interactions occurring on the edge of the pack via suitable nonlinear functions
of the populations in place of the classical bilinear terms coming from the mass action law.
These are therefore Gompertz-like interaction terms, with a fixed exponent, whose value is $\frac 12$. Its value comes from its geometric meaning,
it represents the fact that the perimeter of the patch occupied by the population is one-dimensional, while the patch itself is two-dimensional,
as explained in detail below in Section 2.

The paper is organized as follows. The next Section discusses the model formulation. Section 3
contains the various dynamical systems,
their adimensionalizations, some mathematical preliminaries and the analysis of the equilibria in which one or more populations are absent.
Each of the following Sections investigates instead the coexistence equilibrium respectively for the cases of symbiosis, the two predator-prey
cases and competition. A final discussion concludes the paper.

\section{Modeling pack versus individualistic behavior}\label{sec:2}

The basic ideas underlying modeling herd behavior have been expounded in \cite{1}. Here we recall the main steps for the benefit of the reader.
Consider a population that gathers together. Let $P$
represent its size.
If this population lives on a certain territory of size $A$, the number of individuals staying at the outskirts of the pack is
directly related to the length of
the perimeter of the territory. Therefore its length is proportional to $\sqrt{A}$. Since $P$ is distributed over a two-dimensional
domain, the density square root, i.e. $\sqrt{P}$ will therefore count the individuals at the edge of the territory.

Now let us assume that another population $Q$ intermingles with the one just considered.
We assume that the interactions of the latter occur mainly via the individuals living at the periphery, so that
the interaction term for each individual of the population $Q$ must be proportional to $\sqrt{P}$.
As a result, if $Q$ behaves individualistically, the interactions among the two populations are expressed by $Q\sqrt P$.
Alternatively, if also $Q$ gathers in herds, the interactions will occur at the edge of each bunch of individuals,
and therefore will contain square root terms for both populations. They will thus be modeled via $\sqrt Q \sqrt P$.

Interactions between population can be of different types. They can benefit both, in the case of symbiosis.
Alternatively they can damage both populations, when they compete among themselves directly or for common resources.
Finally, one population receives an advantage from the other one; this happens in the predator-prey situation.

Each of these possible configurations could be subject to pack behavior in one or both populations.
The case of one population gathering in herds while the other one behaves individualistically has already been extensively dealt with in 
\cite{1}. With one exception, that involves pack predation and individual prey, not considered in \cite{1}, we will therefore concentrate
on models involving both populations with individuals sticking together.

\section{Formulation}\label{sec:form}

Let us denote by $P(\tau)$ and $Q(\tau)$ the sizes of two populations in consideration as functions of time $\tau$.
In all the models that follow, the parameters bear the following meaning.
The parameter $r$ is the growth rate of the $Q$ population, with
$K$ or $K_Q$ being its environment's carrying capacity, while, when meaningful, $K_P$ denotes the carrying capacity for the $P$'s.
Further, for the latter $m$ is the natural death rate when the $P$'s are interpreted as predators, models (\ref{mod1}) and (\ref{mod2}),
while it is once more a reproduction rate otherwise, i.e. in the symbiotic (\ref{symb}) and competing (\ref{comp}) cases.
Interaction rates between the two populations are denoted by parameter $q$ for the $Q$ population and by $p$ for the $P$'s.
The following systems will be considered, in which all the parameters are assumed to be nonnegative. 

The symbiotic situation
\begin{eqnarray}\label{symb}
\frac{dQ}{d\tau}=r\left( 1-\frac{Q}{K_Q}\right)Q+q \sqrt{P} \sqrt{Q}, \quad
\frac{dP}{d\tau}=m\left( 1-\frac{P}{K_P}\right)P+p \sqrt{P} \sqrt{Q},
\end{eqnarray}
the predator-prey interactions of pack-individualistic type, for a specialized predator
\begin{equation}\label{mod1}
\frac{dQ}{d \tau}= r\left( 1-\frac{Q}{K}\right) Q-q \sqrt{P}Q, \quad
\frac{dP}{d \tau}= -m P + p \sqrt{P}Q;
\end{equation}
the pack predation - herd defense system, for a specialized predator
\begin{equation}\label{mod2}
\frac{dQ}{d \tau}= r\left( 1-\frac{Q}{K}\right) Q-q \sqrt{P}\sqrt{Q}, \quad
\frac{dP}{d \tau}= -m P + p \sqrt{P}\sqrt{Q};
\end{equation}
and finally the competing case
\begin{eqnarray}\label{comp}
\frac{dQ}{d\tau}=r\left( 1-\frac{Q}{K_Q}\right)Q-q \sqrt{P} \sqrt{Q}, \quad
\frac{dP}{d\tau}=m\left( 1-\frac{P}{K_P}\right)P-p \sqrt{P} \sqrt{Q}.
\end{eqnarray}

Note that when considering predator-prey interactions we need to impose that $p<q$,
since not the whole prey is converted in food for the predators.

\subsection{Models simplification}
As remarked in \cite{1},
singularities could arise in the systems Jacobians when one or both populations vanish.
For the models (\ref{symb}) and (\ref{comp}) we define new variables as follows 
$$
X(t)=\sqrt{\frac{Q(\tau)}{K_Q}}, \ \ Y(t)=\sqrt{\frac{P(\tau)}{K_P}}, \ \ t=\tau\frac{q\sqrt{K_P}}{2\sqrt{K_Q}},
$$
as well as new adimensionalized parameters
$$
a=\frac{K_Q}{K_P} \frac pq, \ \ b=\frac{r\sqrt{K_Q}}{q\sqrt{K_P}}, \ \ c=\frac{m\sqrt{K_Q}}{q\sqrt{K_P}}.
$$  
Therefore the adimensionalized systems read, for (\ref{symb}) 
\begin{equation}\label{adim}
\frac{dX}{dt}=b(1-X^2)X+Y, \quad
\frac{dY}{dt}=c(1-Y^2)Y+aX,
\end{equation}
while for (\ref{comp}) we find
\begin{equation}\label{adimc}
\frac{dX}{dt}=b(1-X^2)X-Y, \quad
\frac{dY}{dt}=c(1-Y^2)Y-aX.
\end{equation}

For the predator-prey cases (\ref{mod1}) and (\ref{mod2}) the substitutions differ slightly. Rescaling 
for
the model (\ref{mod1}) is obtained through
\begin{displaymath}
X=\frac{Q}{K}, \hspace{1cm} Y=\frac{q\sqrt{P}}{m}, \hspace{1cm} t=m\tau ,
\end{displaymath}
and defining the new parameters
\begin{displaymath}
b=\frac{r}{m}, \hspace{1cm} c=\frac{pqK}{2m^2}.
\end{displaymath}
The adimensionalized system can be written as
\begin{equation}\label{mod1s}
\frac{dX}{dt}= b\left( 1-X\right)X - XY, \quad
\frac{dY}{dt}= -\frac{1}{2}Y + cX,
\end{equation}
while in absence of predators, the system reduces just to the first equation.
In this case, easily, the prey follow a logistic growth, toward the adimensionalized carrying capacity $X_1=1$.

For (\ref{mod2}) we have instead
$$
X=\sqrt{\frac{Q}{K}}, \hspace{1cm} Y=\frac{q}{2m}\sqrt{\frac{P}{K}},\hspace{1cm} t=m\tau.
$$
Define now the adimensionalized parameters
$$
e=\frac{r}{2m}, \hspace{1cm} f=\frac{pq}{4m^2}.
$$
The adimensionalized system for $Y>0$ becomes
\begin{equation}\label{mod2s}
\frac{dX}{dt}=e(1-X^2)X-Y, \quad
\frac{dY}{dt}=-\frac{1}{2}Y+fX.
\end{equation}

Note that all the new adimensionalized parameters are combinations of the old nonnegative parameters $r$, $m$, $p$, $q$, $K$;
as a consequence, they must be nonnegative as well.

\subsection{Stability preliminaries}

For the later analysis of the equilibria stability it is imperative to consider the Jacobians of these systems.
We find the following matrices respectively,
for (\ref{adim})
\begin{equation}\label{J_symb}
J^S\equiv
\left(
\begin{array}{cc}
b(1-3X^2) & 1 \\
a & c(1-3Y^2) 
\end{array}
\right)
\end{equation}
and for (\ref{adimc}) 
\begin{equation}\label{J_comp}
J^C\equiv
\left(
\begin{array}{cc}
b(1-3X^2) & -1 \\
-a & c(1-3Y^2) 
\end{array}
\right).
\end{equation}

Considering the predator-prey cases, 
for (\ref{mod1s}) the Jacobian is 
\begin{equation}\label{J_pp1}
J^{PP1}\equiv
\left(
\begin{array}{cc}
b-2bX-Y & -X \\
c & -\frac{1}{2} 
\end{array}
\right),
\end{equation}
while the one for (\ref{mod2s}) reads 
\begin{equation}\label{J_pp2}
J^{PP2}\equiv 
\left(
\begin{array}{cc}
e(1-3X^2) & -1 \\
f & -\frac{1}{2} 
\end{array}
\right) .
\end{equation}

\subsection{The boundary equilibria}\label{sec:b_e}

With this ``geometric'' expression we denote the equilibria in which at least one population vanishes. In fact, they lie on the boundary of the
feasible region of the phase plane, the first quadrant. They need a special care in these kinds of group behavior models, because in
eliminating the singularity we divide by $X$ and $Y$. Therefore all the simplified models (\ref{adim})-(\ref{mod2s})
hold for strictly positive populations. If one population vanishes, no information can gathered by the latter, we rather have to turn
to the original formulations (\ref{symb})-(\ref{comp}).

If one of the two populations disappears the system reduces to one equation.
In this circumstance the surviving population follows a logistic growth toward
its own carrying capacity for the models (\ref{symb}) and (\ref{comp}). The same occurs for the prey in absence of predators in
models (\ref{mod1}) and (\ref{mod2}). In these models when prey are absent, the predators cannot survive.
In fact when $Q=0$ the equation for the predators shows that they exponentially decay to zero.
This makes sense biologically, since these are specialistic predators. Thus in these two models the disappearance of both populations is a possibility.

More generally, the equilibrium corresponding to populations collapse is the origin. Its stability can be analysed by a simple linearization
of the govering equations.

For (\ref{symb}) we find
$$
\frac{dQ}{d\tau}\sim r \sqrt{Q}+q \sqrt{P} >0, \quad
\frac{dP}{d\tau}\sim m \sqrt{P}+p \sqrt{Q}>0.
$$
Thus symbiotic populations cannot both vanish.

For competition, (\ref{comp}), there is a sign change, for which
$$
\frac{dQ}{d\tau}\sim r \sqrt{Q} - q \sqrt{P}, \quad
\frac{dP}{d\tau}\sim m \sqrt{P} - p \sqrt{Q}.
$$
In this case the two populations may disappear, when
\begin{equation}\label{comp_disapp}
\frac mp < \frac {\sqrt{Q}}{\sqrt{P}} < \frac qr.
\end{equation}

For the predator-prey cases, (\ref{mod1}) leads to
$$
\frac{dQ}{d\tau}\sim r \sqrt{Q}>0, \quad
\frac{dP}{d\tau}\sim -m P<0,
$$
so that again the equilibrium is unstable.

In the case (\ref{mod2}) instead we find
$$
\frac{dQ}{d\tau}\sim r \sqrt{Q} - q \sqrt{P}, \quad
\frac{dP}{d\tau}\sim -m \sqrt{P} + p \sqrt{Q}.
$$
Here again both populations under unfavorable circumstances may well disappear, and this happens when
\begin{equation}\label{pp_disapp}
\frac {\sqrt{Q}}{\sqrt{P}} < \min \left\{ \frac mp ,\frac qr \right\}.
\end{equation}

\section{The Symbiotic model}

In view of the preliminary results of Section \ref{sec:b_e}
here as well as in what follows,
we investigate the coexistence of both populations for the simplified models.

\subsection{Analysis}

\begin{figure}[!htb]
\centering
\includegraphics[width=10cm]{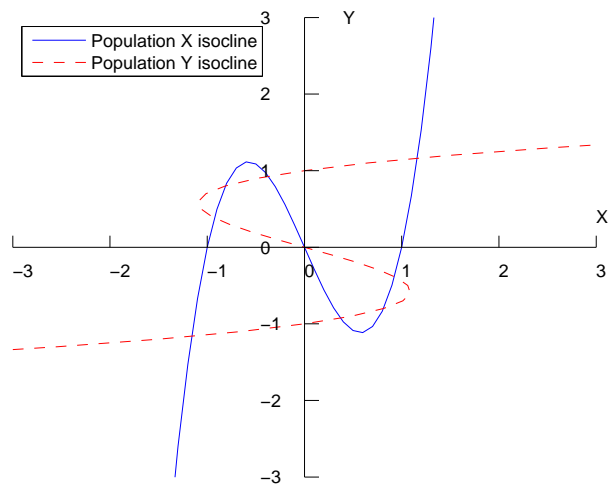}
\caption{Graphical solution of equations system from (\ref{adim}). The $X$ nullcline corresponds to the blue continuous curve $Y=Y(X)$,
conversely The $Y$ nullcline corresponds to the red dashed function $X=X(Y)$.
The phase plane of interest is obviously only the set $\{(X,Y ): X\ge 0, Y \ge 0\}$.
The figure is obtained for the following parameter values $a=0.6$, $b=2.9$, $c=1.7$,
$r=2.9$, $m=1.7$, $p=0.6$, $q=1$, $K_p=10$, $K_q=10$.}
\label{s_phase_whole}
\end{figure}

Looking for the coexistence equilibria we are led to the ninth degree equation
$$
X[a-bc(1-X^2)(1-b^2X^6+2b^2X^4-b^2X^2)]=0
$$
which does not lead to significant results.
However, we then consider a graphical analysis of the system of equations originated by (\ref{adim}).
A typical situation is shown in Figure \ref{s_phase_whole} for an arbitrary choice of the parameter values.

\vspace*{0.25cm}
\begin{proposition}
The internal equilibrium is unique and always feasible.
\end{proposition}

\noindent {\textbf{Proof}}.
As it can be seen, for this parameters choice, all nine roots of the system are real. For other situations, some of the intersections in the second and
fourth quadrant may disappear.
But we are only interested in nonnegative populations $X$ and $Y$.
In view of the behavior of the cubic functions, there will always exist a real intersection between these two functions
in the first quadrant. Moreover, this intersection is unique, leading to the coexistence equilibrium $E_3^S=(X_3^S,Y_3^S)$.

\vspace*{0.25cm}
\begin{proposition}
No Hopf bifurcations can arise at the coexistence equilibrium.
\end{proposition}

\noindent {\textbf{Proof}}.
To have Hopf bifurcations, we need purely imaginary eigenvalues.
This occurs when the trace of the Jacobian vanishes and simultaneously the determinant is positive, i.e.
\begin{equation}\label{RHs}
b(1-3X^2)+c(1-3Y^2)=0, \quad b(1-3X^2)c(1-3Y^2)-a>0.
\end{equation}
It can be easily seen that solving for $b$ from the first condition and substituting into the second one,
we find
$$
a<-c^2(1-3Y^2)^2<0,
$$
which is a contradiction.

\vspace*{0.25cm}
\begin{proposition}
The trajectories of the system (\ref{adim}) are ultimately bounded. $E_3^S$ is globally asympotically stable.
\end{proposition}

\noindent {\textbf{Proof}}. We follow \cite{1} and just outline the proof.
It is enough to take a large enough box $B$ in the first quadrant that contains the coexistence equilibrium.
On the vertical and on the horizontal sides it is easy to show that the dynamical system's flow enters into
the box. The axes cannot be crossed, on biological grounds. Mathematically however, the square root singularity
prevents the right hand side of the dynamical system to be Lipschitz continuous when the corresponding
population vanishes, so that the assumption for the
uniqueness theorem fails on the axes. But as mentioned in the model formulation, we understand that the differential equations
hold only in the interior of the first quadrant, on the coordinate axes they are replaced by corresponding equations in which
the vanishing population is removed.
Thus $B$ is a positively invariant set, from which the first claim follows.
By the Poincar\'e-Bendixson
theorem, since there are no limit cycles, the coexistence equilibrium must be globally asymptotically stable.

\subsection{Comparison with the classical symbiotic case}

The results of the classical case,
\begin{eqnarray}\label{symb_class}
\frac{dQ}{d\tau}=r\left( 1-\frac{Q}{K_Q}\right)Q+q PQ, \quad
\frac{dP}{d\tau}=m\left( 1-\frac{P}{K_P}\right)P+p PQ,
\end{eqnarray}
are summarized in 
\cite{1}. Extensions of classical symbiotic systems have been recently investigated, to models incorporating diseases \cite{NOVA},
or to food chains, \cite{Samrat_NOVA}.
In short, the three equilibria in which at least one population vanishes are unstable,
$\widehat{E_0^S}=(0,0)$, $\widehat{E_1^S}=(K_Q,0)$ and $\widehat{E_2^S}=(0,K_P)$. The coexistence equilibrium 
$$
\widehat{E_3^S}=\left(\frac{K_Qm(r+pK_P)}{rm-pqK_PK_Q},\frac{K_Pr(m+qK_Q)}{rm-pqK_PK_Q}\right)
$$
is unconditionally stable when feasible, i.e. $rm<pqK_PK_Q$.
Note that if $\widehat{E_3^S}$ is infeasible the trajectories are unbounded, which is biologically scarcely possible in view of the
environment's limited resources. 

We now compare the classical model with (\ref{adim}) in order to understand how socialization may boost the mutual benefit of the system's populations.

The symbiotic model (\ref{adim}) has always a stable coexistence equilibrium, while in the classical model $\widehat{E_3^S}$ could be infeasible.

\begin{figure}[!htb]
\begin{minipage}[b]{0.45\linewidth}
\centering
\includegraphics[width=\textwidth]{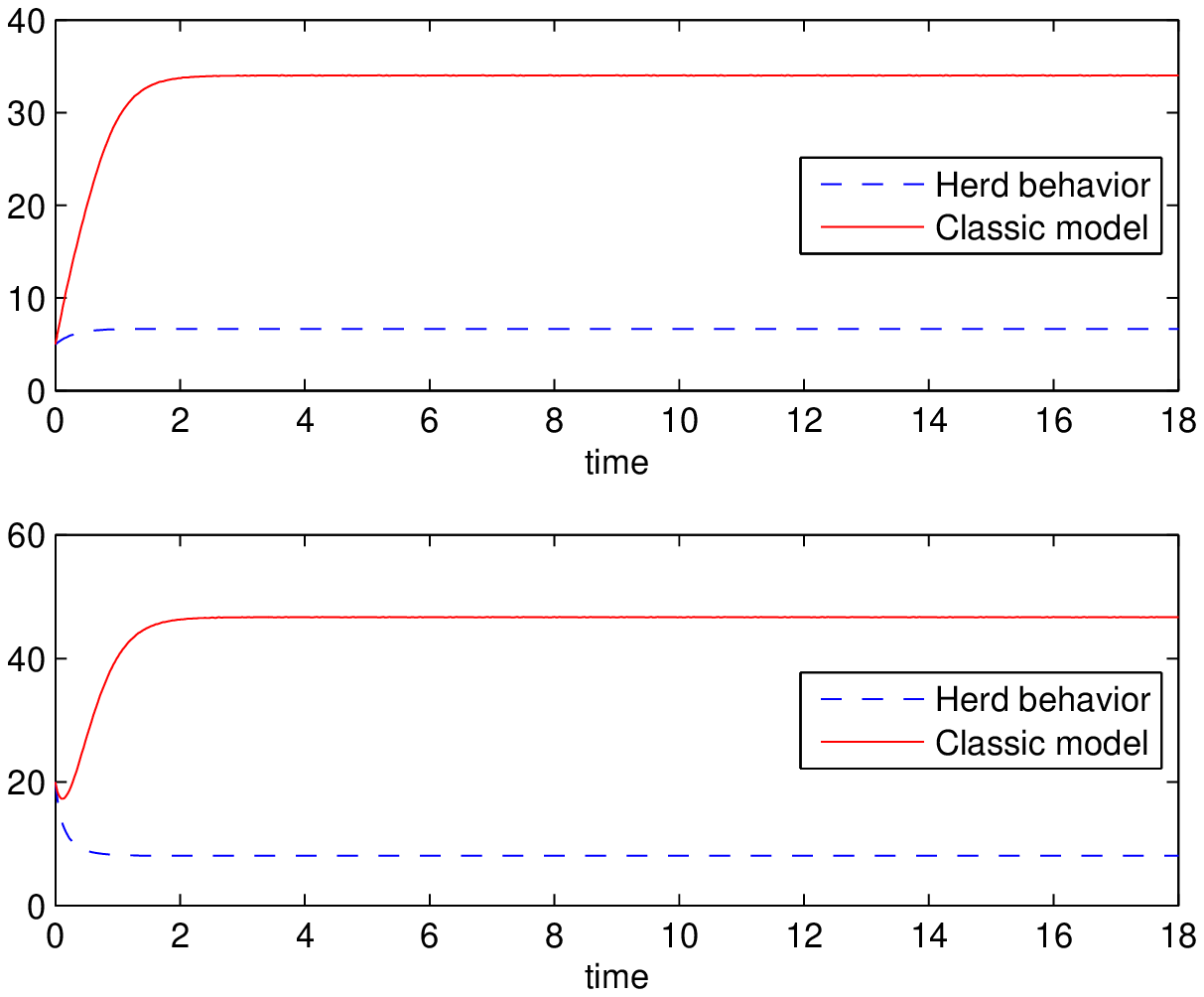}
\end{minipage}
\hspace{0.5cm}
\begin{minipage}[b]{0.45\linewidth}
\centering
\includegraphics[width=\textwidth]{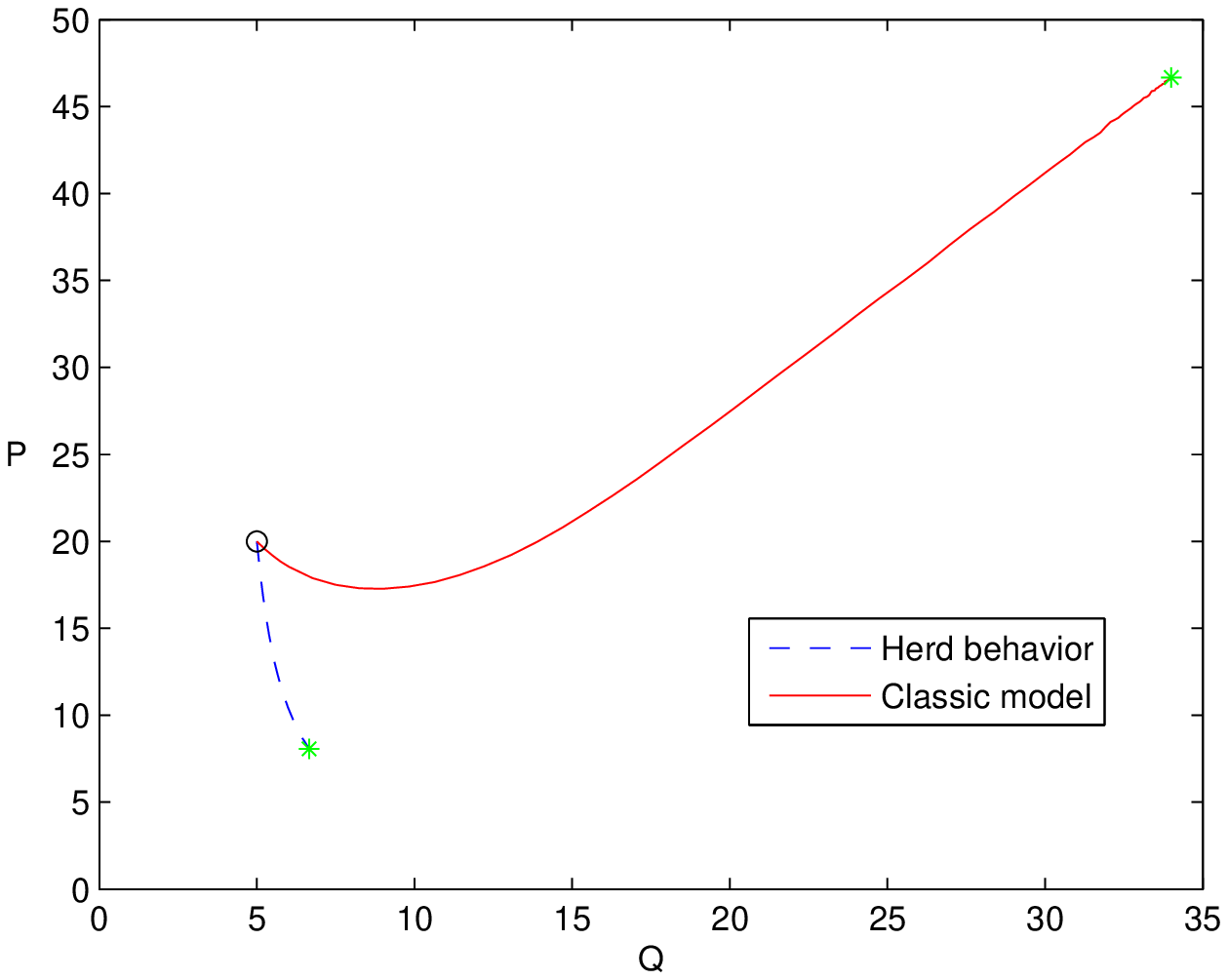}
\end{minipage}
\caption{Left: time series of the systems trajectories; Right: phase plane for classical and new symbiotic model with $r=3$, $m=3$, $K_Q=6$, $K_P=7$, $q=0.3$, and $p=0.5$. Trajectories originate from the same initial condition. The full green dots represent the final equilibrium values.
Rescaled parameter values are $a=1.4286$, $b=9.2582$, $c=9.2582$.}
\label{compare}
\end{figure}

Considering only parameters choices where $\widehat{E_3^S}$ is feasible,
we compare the resulting populations levels for (\ref{adim}) and the classical model.
Taking for both cases $r=3$, $m=3$, $K_Q=6$, $K_P=7$, $q=0.3$, and $p=0.5$, the behaviors are shown in Figure \ref{compare}. Starting from
the same initial conditions, different equilibria are reached. 

Clearly the population level is higher in the classical model. The numerical values we obtained are $Q=6.66$, $P=8.06$
for the herd model and $Q=33.99$, $P=46.69$ for the classical model.

This makes sense, since in symbiotic models the benefit comes from the mutual interactions between populations.
If the latter are scattered in the environment it is more likely for each individual of one population to get in contact with one of the other.
On the other hand, when herd behaviour is exhibited, only individuals on the outskirts interact with the other population
and as a consequence the innermost individuals receive less benefit since they hardly have the chance to meet the other population.

\section{The predator-prey cases}

Here, we let $P$ represent the density of the predators and $Q$ denote the prey population.

\subsection{Pack predation and individualistic prey behavior}

We consider now (\ref{mod1}) in the adimensionalized form (\ref{mod1s}).
We can immediately show that the trajectories are bounded.

\vspace*{0.25cm}
\begin{proposition}
All populations in (\ref{mod1s}) are bounded.
\end{proposition}

\noindent {\textbf{Proof}}. 
Introducing the environment total population, $Z(t)=X(t)+Y(t)$ and summing the equations in (\ref{mod1s}), we have
\[
\frac{dZ}{dt}=-\frac{1}{2}Y+cX+bX-bX^2-XY =-\frac{1}{2}Z+\left(c+b+\frac{1}{2}-bX-Y\right)X.
\]
Take the maximum of the parabola in $X$ on the right hand side, to obtain
\[
\frac{dZ}{dt}+\frac{1}{2}Z \leq\left(c+b+\frac{1}{2}-bX\right)X  \leq \frac{\left(c+b+\frac{1}{2}\right)^2}{4b} \equiv \bar{M}.
\]
The above differential inequality leads to
$$
Z(t) \leq e^{-\frac{1}{2}t}+2\bar{M}\left(1-e^{-\frac{1}{2}t}\right)\leq 1+2\bar{M}=M.
$$
Because the total population is bounded, also each individual population $X$ and $Y$ is bounded as well.

\vspace*{0.25cm}

Here 
the coexistence equilibrium is always feasible,
$$
E_2=\left(\frac{b}{b+2c},\frac{2bc}{b+2c}\right).
$$

\vspace*{0.25cm}
\begin{proposition}
The coexistence equilibrium $E_2$ is always locally asymptotically stable.
\end{proposition}

\begin{figure}[!htb]
\centering
\includegraphics[width=8cm]{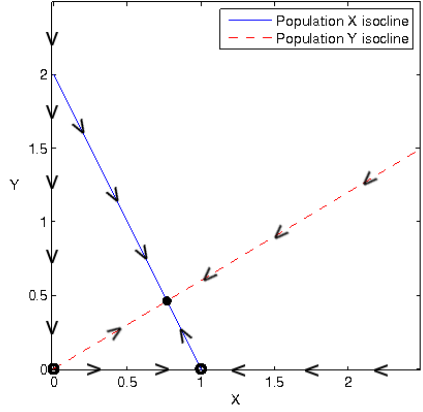}
\caption{Phase plane sketch of the model (\ref{mod1}) with parameters values
$b=2$, $c=0.3$, $r=0.6$, $m=0.3$, $p=0.0072$, $q=1.5$, $K=5$.}
\label{phase1}
\end{figure}

\noindent {\textbf{Proof}}.
If $J^{PP1}_2$ denotes the Jacobian matrix (\ref{J_pp1}) evaluated at $E_2$,
the Routh-Hurwitz criterion gives
\[
\det(J^{PP1}_2)=-\frac{1}{2}b+\frac{b^2+2bc}{b+2c}=\frac{1}{2}b>0, \quad
\textrm{tr}(J^{PP1}_2)=-\frac{1}{2}+b-\frac{2b^2+2bc}{b+2c}=-\frac{2b^2+2c+b}{2(b+2c)}<0.
\]
Both conditions hold so that the eigenvalues have negative real part and $E_2$ is always a stable equilibrium.
The phase plane picture also supports this conclusion as well, Figure \ref{phase1}.

\vspace*{0.25cm}
\begin{proposition}
The coexistence equilibrium $E_2$ is also globally asymptotically stable.
\end{proposition}

\noindent {\textbf{Proof}}. 
It follows the outline of Proposition 3. More formally,
consider the point $N^*\equiv (N,2cN)$, with $N>1$, on the isocline through the origin in Figure \ref{phase1}. The compact set
$\Omega$, identified by the rectangle having $N^*$ and the origin as opposite vertices, is positively invariant.
On its right vertical side indeed we have
$$
\frac {dX}{dt}\vert _{X=N}= b(1-N)N-2cN^2<0
$$
and the system's trajectories must enter into $\Omega$ from the right. On the upper side, for $X<N$, we have
$$
\frac {dY}{dt}\vert _{Y=2cN}= -\frac 12 2cN+cX <0.
$$
Here the trajectories of (\ref{mod1s}) enter into $\Omega$ from above. By the Poincar\'e-Bendixson theorem, global stability follows.

\vspace*{0.25cm}

Note that Hopf bifurcations cannot arise here, since $\textrm{tr}(J^{PP1}_2)<0$ is a strict inequality.

\subsection{Prey herd behavior}

We focus now on (\ref{mod2s}).

\vspace*{0.25cm}
\begin{proposition}
All populations in (\ref{mod2s}) are bounded.
\end{proposition}

\noindent {\textbf{Proof}}. The steps are the same as for Proposition 4, with minor changes. The differential inequality here becomes
\[
\frac{dZ}{dt}+\frac{3}{2}Z \leq \left(\frac{2}{3}e+\frac{2}{3}f+1\right)\sqrt{\frac{2e+2f+3}{6e}} \equiv \bar{M}.
\]
where the last estimate follows by taking the maximum of the cubic in $X$.
We then find
$$
Z(t) \leq e^{-\frac{3}{2}t}+\frac{2}{3}\bar{M}\left(1-e^{-\frac{3}{2}t}\right)\leq 1+\frac{2}{3}\bar{M}=M,
$$
providing a bound on both populations as well as on each subpopulation.

\vspace*{0.25cm}

\begin{figure}[!htb]
\begin{minipage}[b]{0.45\linewidth}
\centering
\includegraphics[width=\textwidth]{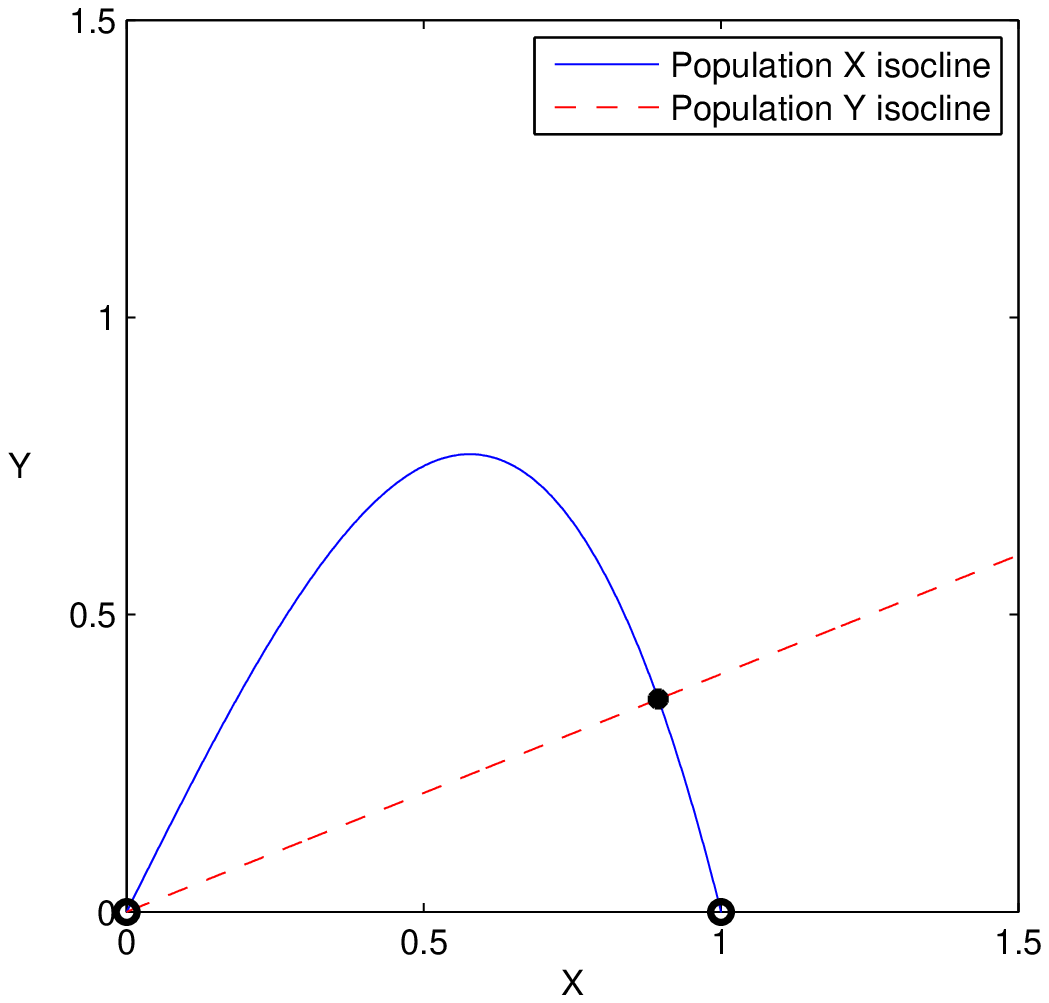}
\caption{Phase plane of model (\ref{mod2s}) with $e\geq2f$, both $\widehat E_0$ and $\widehat E_2$ exist.
Parameter values: $e=2$, $f=0.2$, $r=0.5$, $m=0.125$, $p=0.5$, $q=0.025$, $K=10$.}
\label{phase2}
\end{minipage}
\hspace{0.5cm}
\begin{minipage}[b]{0.45\linewidth}
\centering
\includegraphics[width=\textwidth]{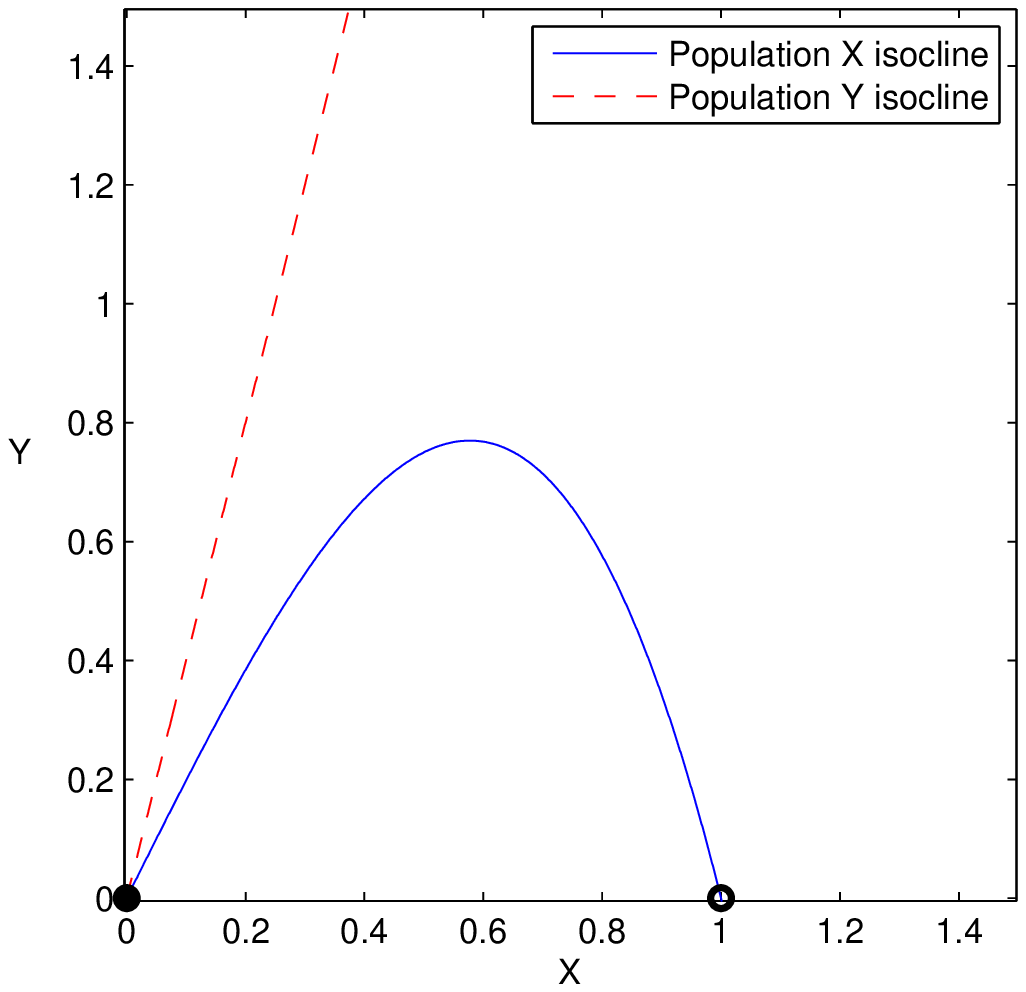}
\caption{Phase plane of model (\ref{mod2s}) with $e<2f$, $\widehat E_2$ is infeasible.
Parameter values: $e=2$, $f=2.0$, $r=0.5$, $m=0.125$, $p=0.5$, $q=0.25$, $K=10$.}
\label{phase22}
\end{minipage}
\end{figure}

The coexistence equilibrium
$$
\widehat E_2=\left(\sqrt{\frac{e-2f}{e}},2f\sqrt{\frac{e-2f}{e}}\right),
$$
is feasible for
\begin{equation}\label{2E2_feas}
e\geq2f.
\end{equation}
Figures \ref{phase2} and \ref{phase22} illustrate geometrically the two situations in which $\widehat E_2$ is feasible and when it is infeasible.

Recalling that in the case of (\ref{mod2}) the origin might be stable, (\ref{pp_disapp}),
and that when (\ref{2E2_feas}) becomes an equality the coexistence equilibrium vanishes, we have the following result.

\vspace*{0.25cm}
\begin{proposition}
There is a transcritical bifurcation
for which $\widehat E_2$ emanates from the origin $\widehat E_0$ when the parameter $e$ raises up to attain the critical value $e^*=2f$.
\end{proposition}

\noindent {\textbf{Proof}}.
The characteristic polynomial at the origin $\widehat E_0$ is
$$
\lambda^2+\left(\frac{1}{2}-e\right)\lambda+f-\frac{1}{2}e=0.
$$
The Routh-Hurwitz stability conditions then become
\begin{equation}\label{2E0_stab}
2f>e, \hspace{1cm} e<\frac{1}{2}.
\end{equation}
The second claim follows comparing the first inequality in (\ref{2E0_stab}) with (\ref{2E2_feas}). In fact, at $e^*$ the origin becomes
unstable, while instead $\widehat E_2$ becomes feasible.

\vspace*{0.25cm}
\begin{proposition}
Coexistence for the system (\ref{mod2s}) is a locally asymptotically stable equilibrium either if
$2f<e<\frac{1}{2}$ and (\ref{tr}) holds; or if $e>\max \left\{\frac{1}{2}, 3f-\frac{1}{4}\right\}$ and (\ref{tr}) holds.
But if $\frac{1}{2}<e<3f-\frac{1}{4}$ we find that (\ref{tr}) is not true and $E_2$ is unstable.
\end{proposition}

\begin{figure}[!htb]
\centering
\includegraphics[width=8cm]{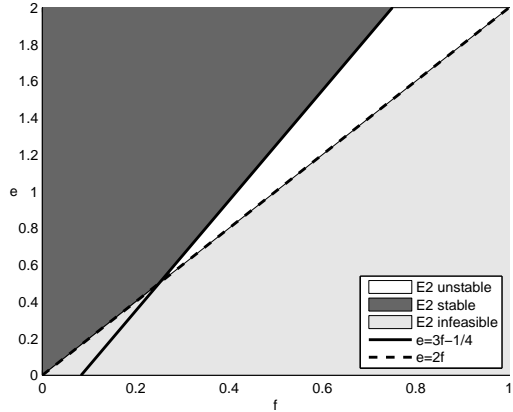}
\caption{Region of the $f-e$ parameter space in which the coexistence equilibrium of (\ref{mod2s}) is stable.}
\label{fig-stab}
\end{figure}

\noindent {\textbf{Proof}}.
Let the Jacobian evaluated at $\widehat E_2$ be denoted by $J^{PP2}_2$.
The Routh-Hurwitz conditions are now $\det(J^{PP2}_2)=e-2f>0$, which always holds by the feasibility condition (\ref{2E2_feas}), and
\begin{equation}\label{tr}
{\textrm{tr}}(J^{PP2}_2)=-2e+6f-\frac{1}{2}<0. 
\end{equation}
There are a few different situations for (\ref{tr}), represented in Figure \ref{fig-stab}.

\vspace*{0.25cm}
\begin{proposition}
When locally asymptotically stable, the equilibrium $\widehat E_0$
is also globally asymptotically stable.
\end{proposition}

\noindent {\textbf{Proof}}. Proceeding as for (\ref{mod1s}), we take the point $\widehat L_*\equiv (L,2fL)$ with 
$$
L> \max \left\{ 1, \frac e{3\sqrt 3 f} \right\}.
$$
The rectangle $\widehat \Omega$ with the origin and $\widehat L_*$ as opposite vertices is a positively invariant set.
Recall that the part of the horizontal sides of interest here is $X<L$.
On the right vertical
and upper horizontal sides of $\widehat \Omega$ indeed, we have
$$
\left. \frac {dX}{dt}\right| _{X=L}=e(1-L^2)L-Y<0, \quad
\left. \frac {dY}{dt}\right|_{Y=2fL}=-f(L-X)<0.
$$
All trajectories thus enter into $\widehat \Omega$.
The only locally asymptotically stable equilibrium in its interior
must also be globally asymptotically stable by the Poincar\'e-Bendixson theorem.

\vspace*{0.25cm}

We summarize the equilibria of system (\ref{mod2s}) in the following table.
\begin{center}
\begin{tabular}{lccc}
\hline
Parameter conditions & $\widehat E_0$ & $\widehat E_2$ & Bifurcation \\
\hline
\vspace{0.2cm}

$e<\frac{1}{2}$ \hspace{0.3cm} $f>\frac{e}{2}$ & stable & \hspace{0.3cm} unfeasible \\ \vspace{0.2cm}

$e<\frac{1}{2}$ \hspace{0.3cm} $e^*=ef$ & & & Transcritical  \\ \vspace{0.2cm}

$e<\frac{1}{2}$ \hspace{0.3cm} $f<\frac{e}{2}$ & unstable & \hspace{0.3cm} stable   \\ \vspace{0.2cm}

$e>\frac{1}{2}$ \hspace{0.3cm} $e>3f-\frac{1}{4}$ & unstable & \hspace{0.3cm} stable  \\ \vspace{0.2cm}

$e>\frac{1}{2}$ \hspace{0.3cm} $e=e^{\dagger}=3f-\frac{1}{4}$ \hspace{0.3cm} & & & Hopf \\ \vspace{0.2cm}

$e>\frac{1}{2}$ \hspace{0.3cm} $2f<e<3f-\frac{1}{4}$ \hspace{0.3cm} & unstable & \hspace{0.3cm} unstable \\ \vspace{0.2cm}

$e>\frac{1}{2}$ \hspace{0.3cm} $f>\frac{e}{2}$ & unstable & \hspace{0.3cm} unfeasible  \\
\hline
\end{tabular}
\end{center}

\vspace*{0.25cm}
\begin{proposition}
The system (\ref{mod2s}) admits a Hopf bifurcation at the coexistence equilibrium when the bifurcation parameter $e$ crosses the critical value
\begin{equation}\label{hpf}
e^{\dagger}=3f-\frac{1}{4}.
\end{equation}
\end{proposition}

\noindent {\textbf{Proof}}.
In addition to the transcritical of Proposition 8,
we show now that special parameters combinations originate Hopf bifurcations near $\widehat E_2$.
Recall that purely imaginary eigenvalues are needed, and this occurs when the trace of the Jacobian vanishes. Thus (\ref{tr}) must
become an equality and the constant term is positive, $\det(J^{PP2}_2)=e-2f>0$. But the latter holds from (\ref{2E2_feas}).

\begin{figure}[!htb]
\centering
\includegraphics[width=7cm]{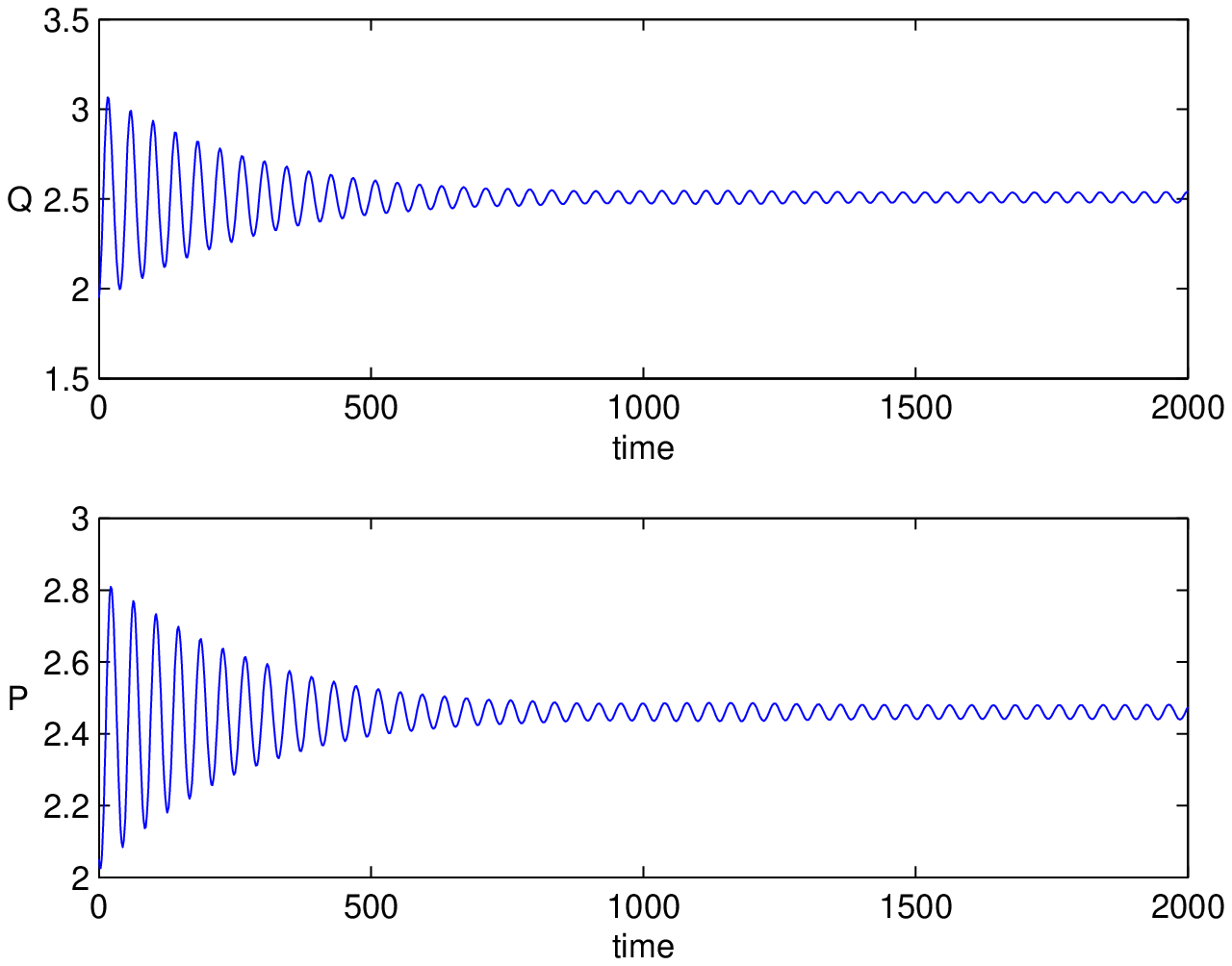}
\includegraphics[width=7cm]{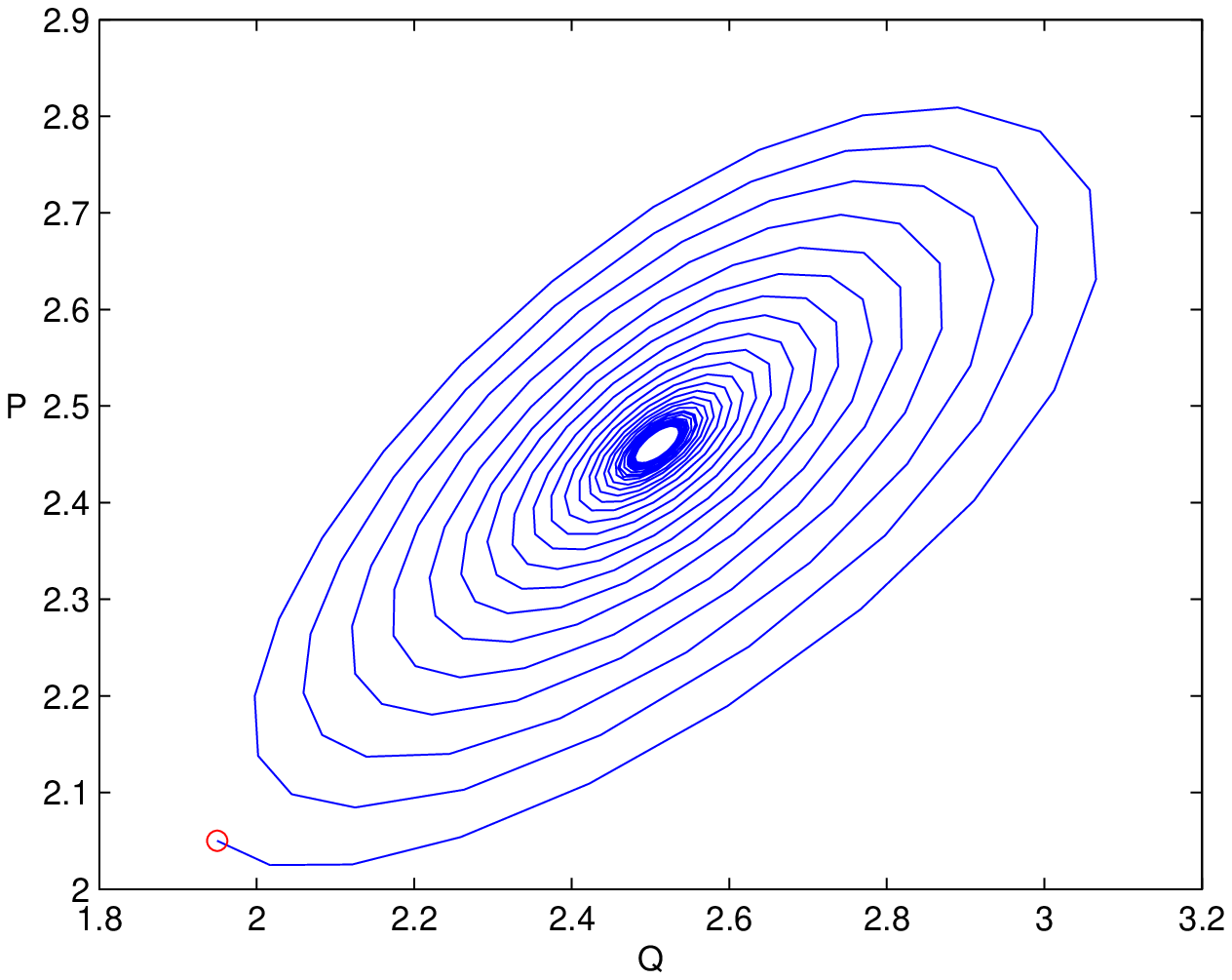}
\caption{Left: time series of the system trajectories (\ref{mod2}); Right: corresponding limit cycle in the phase plane.
The original parameter values are 
$r=0.76$, $m=0.2999$, $p=0.297$, $q=0.607$, $K=12$,
with coexistence equilibrium
$E_2=(2.5085,2.4602)$;
they correspond to $e=1.2671$, $f=0.5011$
in the rescaled model (\ref{mod2s}).}
\label{fig-cycles}
\end{figure}

\vspace*{0.25cm}
This result is observed in Figure \ref{fig-stab}, where the thick straight line indicates the critical parameter values.
Figure \ref{fig-cycles} shows the limit cycles for the dimensionalized model (\ref{mod2}), letting the simulation run for long times to show
that the oscillations are indeed persistent.

\vspace*{0.25cm}
\begin{proposition}
The system (\ref{mod2}) admits trajectories for which the prey go to extinction in finite time, if their initial conditions
lie in the set
\begin{equation}\label{extinction_set}
\Xi=\{ (Q,P):Q>0, P>Q q^{-2} (m+r)^2\}.
\end{equation}
\end{proposition}

\noindent {\textbf{Proof}}.
We follow with suitable modifications the argument exposed in \cite{PV}.
From the second equation in (\ref{mod2}) we get the differential inequality
\begin{equation}\label{d_ineqP}
\frac {dP}{d\tau} \ge -mP
\end{equation}
from which
$P(\tau)\ge \widehat P(\tau) = P_0 \exp(-m\tau)$, where the latter function denotes the solution of the differential equation corresponding to (\ref{d_ineqP}),
with $\widehat P(0)=P(0)$. From the first equation in (\ref{mod2}) we have further
\begin{equation}\label{d_ineqQ}
\frac {dQ}{d\tau} \le rQ-q\sqrt {P}\sqrt {Q}\le rQ-q\sqrt {\widehat P}\sqrt {Q}.
\end{equation}
Let $\widehat Q(\tau)$ denote the solution of the differential equation obtained from (\ref{d_ineqQ}) using the rightmost term,
with $\widehat Q(0)=Q(0)$. It follows that
$Q(\tau)\le \widehat Q(\tau)$. Using the integrating factor $W(\tau)=\widehat Q(\tau) \exp(-r\tau)$,
we obtain
$$
\sqrt {W(\tau)}=\sqrt {W(0)} -\frac {q\sqrt {P(0)}}{m+r}  \left[ 1-\exp\left(-\frac {m+r}2 \tau \right)\right].
$$
The last term on the right is an increasing function of $\tau$, so that there is a $\tau^*$ for which $W(\tau^*)=\widehat Q(\tau^*)=0$ if and only if
\begin{equation}\label{c_init}
\sqrt {W(0)} < \frac {q\sqrt {P(0)}}{m+r}.
\end{equation}
Since $W(0) =Q(0)$, we have $\widehat Q(\tau^*)=0$ if the following inequality for the initial conditions of the trajectories is satisfied,
$$
\sqrt {P(0)} > \frac {m+r}q \sqrt {W(0)},
$$
from which the set $\Xi$ given in (\ref{d_ineqP}) is immediately obtained.

\subsection{Comparison with the classical predator-prey model}

In order to compare these results
quantitatively, we consider also the classical model with logistic correction.
If we rescale it, however, since it does not contain
the square root terms, we would find a different adimensionalization, rendering the comparison difficult.
Thus we rather return to the original formulations
(\ref{mod1}) and (\ref{mod2}).

The classical Lotka-Volterra model with logistic correction for the prey,
\begin{eqnarray}\label{mod1_class}
\frac{dQ}{d\tau}=r\left( 1-\frac{Q}{K_Q}\right)Q-q PQ, \quad
\frac{dP}{d\tau}=-mP+p PQ,
\end{eqnarray}
has two attainable equilibria:
the prey-only equilibrium $(K,0)$, from which via a transcritical bifurcation at $m^{\dagger}=pK$
coexistence $C_*$ arises, feasible for $pK>m$ and always stable.
Its explicit representation follows, together with the one of coexistence for the system with individualistic hunting and prey herd response, \cite{1},
in dimensionalized form,
$$
C_*\equiv \left( \frac mp, \frac rq\left( 1-\frac m{pK}\right)\right), \quad
\widetilde E_2=\left(\frac{m^2}{p^2},\frac{mr}{pq}\left( 1-\frac{m^2}{p^2K}\right)\right).
$$
The dimensional form of the coexistence equilibria of the two models presented here, respectively for (\ref{mod1}) and (\ref{mod2}) are
$$
E_2\equiv \left( K\frac {mr}{mr+pqK}, \frac {p^2}{m^2} K^2\left(  \frac  {mr}{mr+pqK}\right)^2\right), \quad
\widehat E_2\equiv \left( K\left( 1-\frac {pq}{mr}\right), \frac {p^2}{m^2}K\left( 1-\frac {pq}{mr}\right)\right).
$$

The equilibrium prey populations of the first two models depend only on the system parameters $m$ and $p$, i.e. the predators' mortality and predation efficiency.
Thus they are independent of their own reproductive capabilities and of the environment carrying capacity.
Further, when the predators' hunting efficiency is larger than the predators' own mortality, i.e. $m<p$,
the equilibrium prey value is  much lower if they gather in herds, i.e. in $\widetilde E_2$, while on the contrary the predators attain instead higher values,
again at $\widetilde E_2$.
Conversely, when $m>p$ the prey grouping together, $\widetilde E_2$,
allows higher equilibrium numbers than for their individualistic behavior; the predators instead
settle at lower values if the prey use a defensive strategy, $\widetilde E_2$, and higher ones with individualistic prey behavior, at $C_*$.

For (\ref{mod1}) and (\ref{mod2}), i.e. with coordinated hunting, these values involve also the prey own intrinsic characteristics.
In particular for (\ref{mod2})
the ratio of the predators' hunting efficiency $p$ versus their mortality $m$ determines if the predators at equilibrium
will be more than the prey, see $\widehat E_2$. 

A similar result possibly extends  for the model of pack hunting coupled with loose prey, (\ref{mod1}),
but at $E_2$ the predators population at equilibrium contains the prey population squared
and in principle the latter may not exceed $1$, so that the conclusion
would not be immediate.
Indeed,
at the equilibria $E_2$ and $\widehat E_2$, the
prey populations are the multiplication of the fractions in the brackets, always smaller than $1$, by the carrying capacity $K$, which
may or not be large.
The result could indeed give a population
smaller than $1$. This in principle is not a contradiction, because the population need not necessarily be counted by
individuals, but rather its size could be measured by the weight of its biomass.

\section{The competition model}
We now deal with (\ref{comp}) in the rescaled version (\ref{adimc}).

\subsection{Analysis}

\begin{figure}[!htb]
\centering
\includegraphics[width=0.48\textwidth]{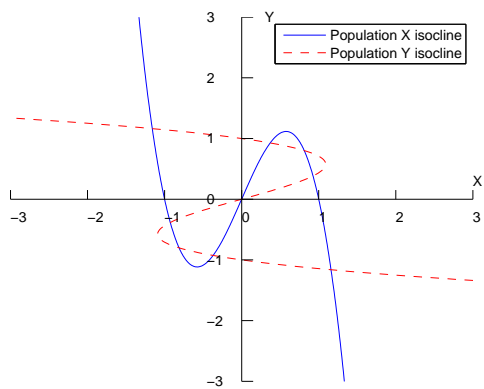}
\caption{Graphical solution of equations system from (\ref{adimc}) for the functions $Y_{[1]}(X)$ and $X_{[2]}(Y)$.
Parameter values: $a=0.6$, $b=2.9$, $c=1.7$, $r=2.9$, $m=1.7$, $p=0.6$, $q=1$, $K_p=10$, $K_q=10$.}
\label{c_phase_whole}
\end{figure}

The coexistence equilibria are the roots of the eighth degree equation
\begin{equation}\label{8th_deg_c}
cb^3X^8-3cb^3X^6+3cb^3X^4-cb(b^2+1)X^2-a+cb=0.
\end{equation}
A better interpretation treats the problem as an intersection of cubic functions,
\begin{equation}\label{cubic_c}
Y_{[1]}(X)=b(1-X^2)X, \quad X_{[2]}(Y)=\frac{c}{a}(1-Y^2)Y.
\end{equation}
Depending on the behavior of the cubic functions, there could be either three intersections (the origin and one each in the second and fourth quadrants)
or five (the previous ones and one more in the first and third quadrants), or nine.
The latter configuration is graphically shown in Figure \ref{c_phase_whole}.
The feasible coexistence equilibria are just the intersections in  the first quadrant.
Note that no intersections in the first quadrant exist when the slopes at the origin of the two cubic functions (\ref{cubic_c}) satisfy the inequality
$Y_{[1]}'(0)<Y_{[2]}'(0)$, the latter denoting of course the inverse function of $X_{[2]}(Y)$. This condition, rephrased in terms of the parameters, becomes
\begin{equation}\label{a_bc}
a>bc.
\end{equation}

\begin{figure}[!htb]
\begin{minipage}[b]{0.3\textwidth}
\includegraphics[width=1.1\textwidth]{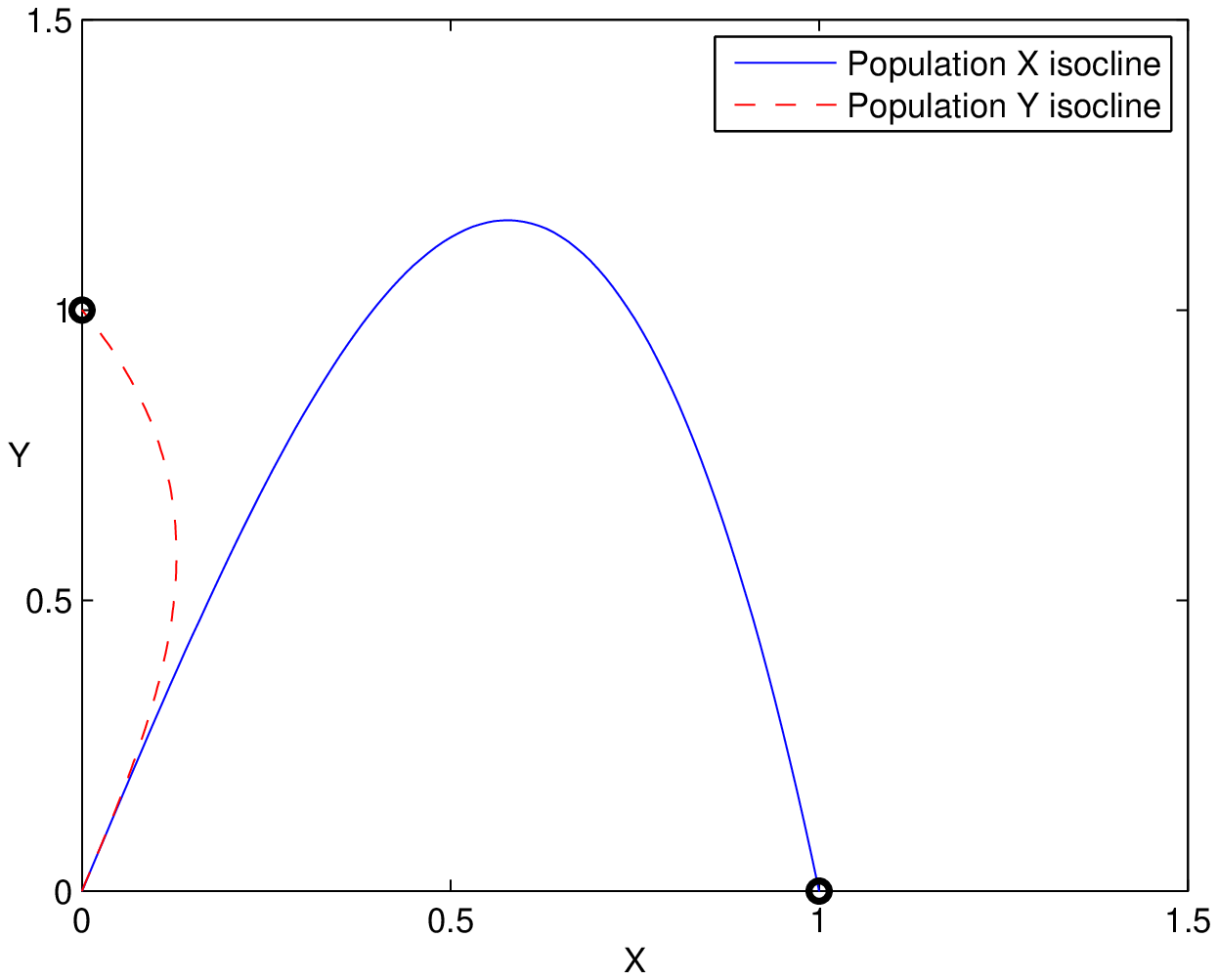}
\label{0equilibrium}
\end{minipage}
\begin{minipage}[b]{0.3\textwidth}
\includegraphics[width=1.1\textwidth]{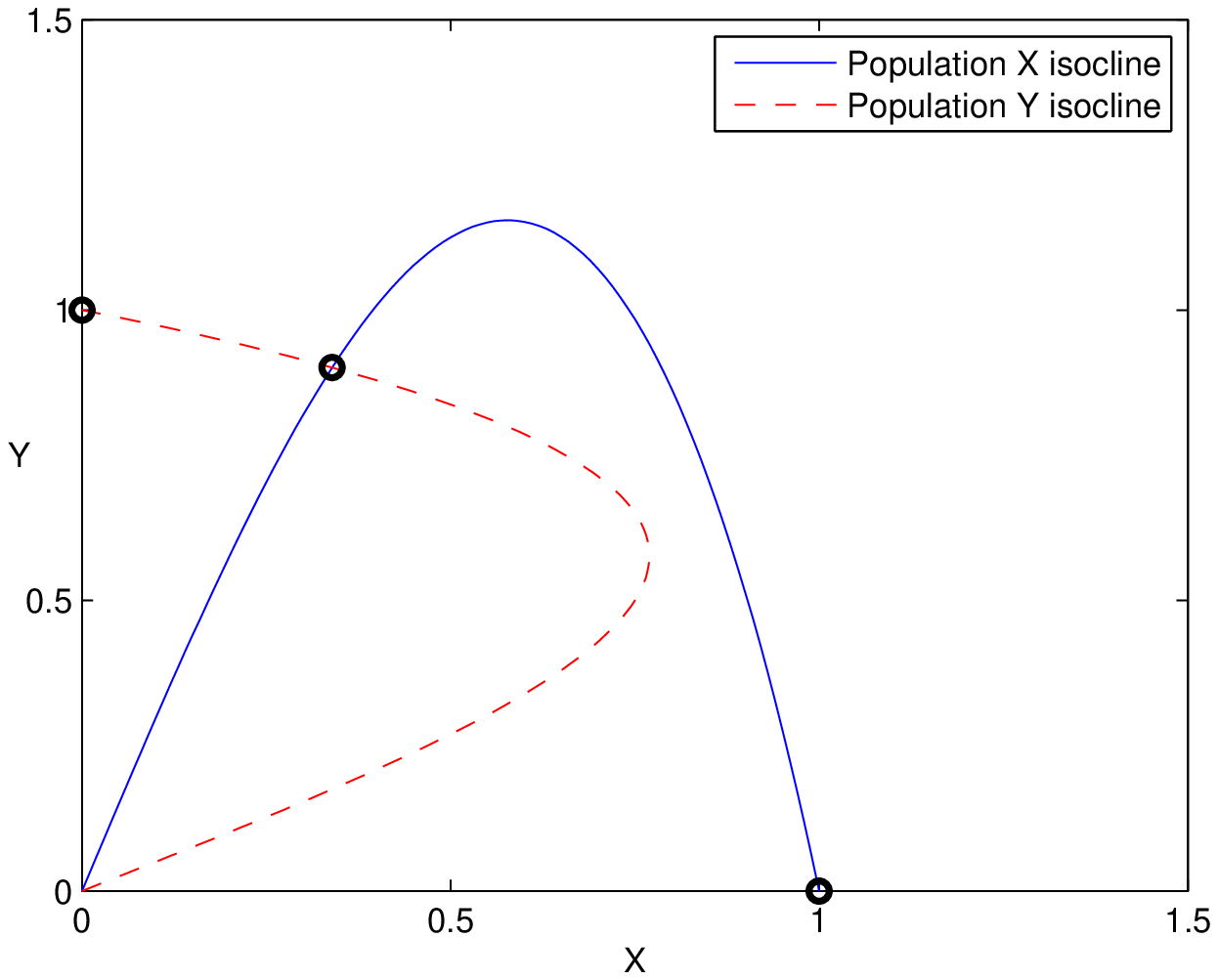}
\label{1equilibrium}
\end{minipage}
\begin{minipage}[b]{0.3\textwidth}
\includegraphics[width=1.1\textwidth]{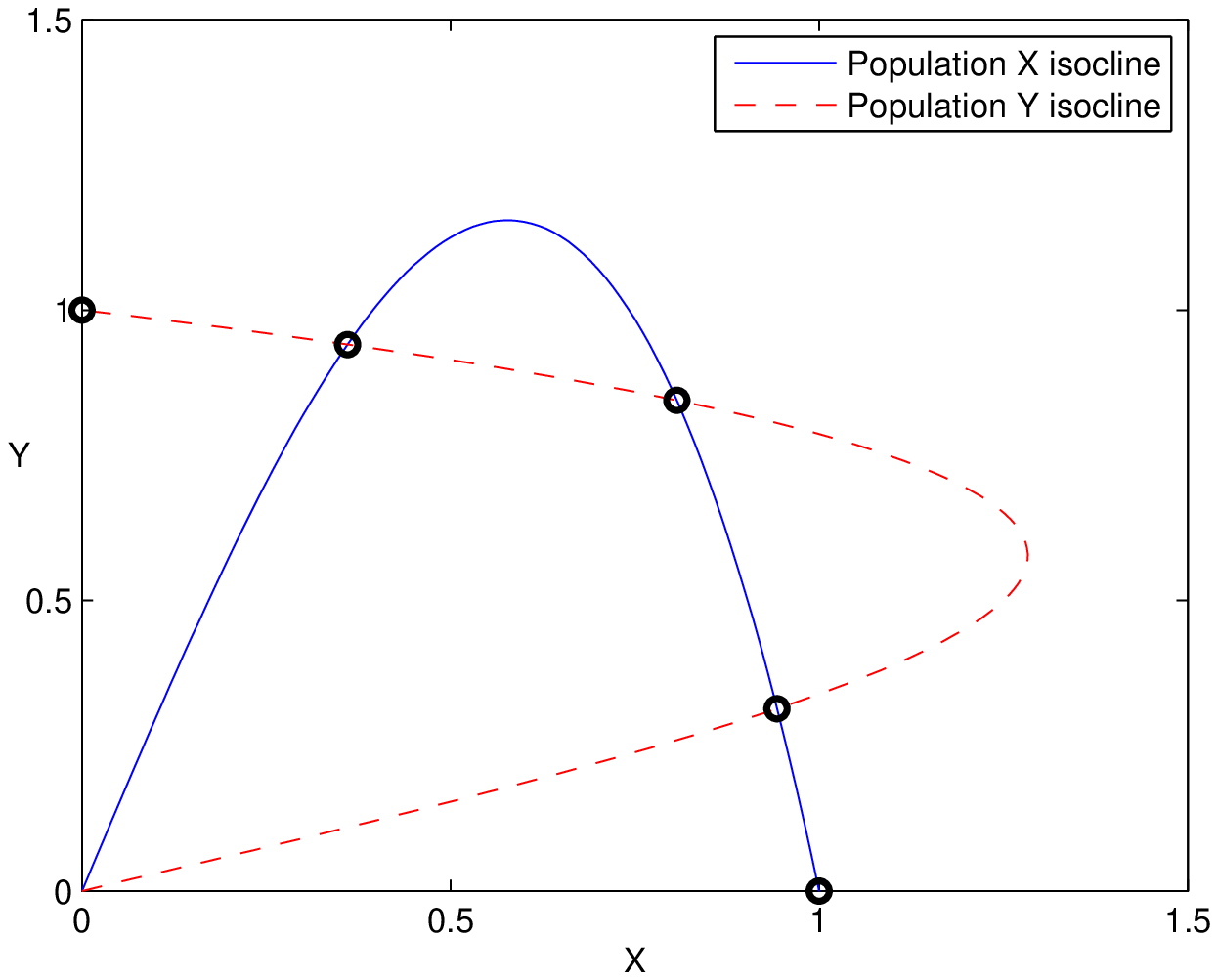}
\label{3equilibria}
\end{minipage}
\caption{Referring back to Figure \ref{c_phase_whole}, we show here the
coexistence equilibria possible scenarios: left $a>bc$ no feasible equilibrium exists for the parameter values
$a=3$, $b=3$, $c=1$, $r=0.9$, $m=0.3$, $p=0.9$, $q=0.3$, $K_p=10$, $K_q=10$; center $b<\frac{3\sqrt{3}}{2}$ and $c<\frac{3\sqrt{3}}{2}a$,
just one feasible equilibrium $E_3^C$, for the parameter values
$a=3$, $b=3$, $c=6$, $r=0.9$, $m=1.8$, $p=0.9$, $q=0.3$, $K_p=10$, $K_q=10$;
right $b>\frac{3\sqrt{3}}{2}$ and $c>\frac{3\sqrt{3}}{2}a$, for the parameter values
$a=3$, $b=3$, $c=10$, $r=0.9$, $m=3$, $p=0.9$, $q=0.3$, $K_p=10$, $K_q=10$. The three equilibria $E_4^C$, $E_3^C$ and $E_5^C$ are ordered
left to right, for increasing values of their abscissae.}
\label{scenario}
\end{figure}

Thus, for $a>cb$ there is at most one real positive root, the one corresponding to the intersection in the fourth quadrant, that
is however not feasible, see the left frame in Figure \ref{scenario}.
Thus no coexistence equilibrium arises.

To better analyse the situation, we apply Descartes' rule of signs to (\ref{8th_deg_c}). There are three sign variations, since the first
four coefficients have alternating signs. The last one must be positive, because having already ruled out the case (\ref{a_bc}), we are left with
$a<cb$. Descartes' rule shows that in this case
there are at most 4 real positive roots. Recall that these roots correspond to the abscissae of the intersections of
the curves (\ref{cubic_c}). As discussed above we know
that one positive root corresponds to the intersection that always exists in the fourth quadrant, Figure \ref{c_phase_whole}. This root must then
be excluded.
As a consequence in this case we have just one or three coexistence equilibria, see the center and right frames in Figure \ref{scenario}.

Sufficient conditions for three versus one equilibria to exist is that the cubic functions (\ref{cubic_c})
have maximum $Y$-coordinate and $X$-coordinate respectively in the first quadrant greater than 1. This happens when
both the following conditions hold
$$
b>\frac{3\sqrt{3}}{2}, \quad c>\frac{3\sqrt{3}}{2}a.
$$

The three possible situations are shown in Figure \ref{scenario}.

In summary we have the following result.

\vspace*{0.25cm}
\begin{proposition}
If $a>bc$ no feasible coexistence equilibria exist. If $a<bc$ at least one feasible
equilibrium exists, $E_3^C=(X_3^C,Y_3^C)$. Further, in such case,
$b>\frac{3\sqrt{3}}{2}$ and $c>\frac{3\sqrt{3}}{2}a$ are sufficient conditions
for three equilibria to exist, i.e. $E_4^C$, $E_3^C$ and $E_5^C$,
ordered for increasing values of their abscissae.
\end{proposition}

\vspace*{0.25cm}
\begin{proposition}
The trajectories of the system (\ref{adimc}) are ultimately bounded.
\end{proposition}

\noindent {\textbf{Proof}}. Observe that $X$ decreases when $Y\le bX(1-X^2)$ and similarly $Y$ decreases for $X\le ca^{-1}Y(1-Y^2)$. This in the phase plane
corresponds to having the flow entering a suitable box $\Omega^C$ with one corner in the origin and the opposite one $\Omega^C_B=(X_B,Y_B)$
of size large enough to contain the vertices of the cubics in all
cases of Figure \ref{scenario}. Thus we can take $X_B\ge \max\{ 1, X_V\}$, $Y_B\ge \max\{ 1, Y_V\}$, where $X_V$ and $Y_V$ denote respectively
the relative maxima heights of the cubics.

\vspace*{0.25cm}
\begin{proposition}
The equilibria for which either one of the conditions hold
\begin{equation}\label{p13_assump}
X<\frac {\sqrt 3}3, \quad Y<\frac {\sqrt 3}3
\end{equation}
are unstable.
\end{proposition}

\vspace{0.5cm}
\noindent {\textbf{Proof}}.
If both (\ref{p13_assump}) hold,
the first Routh-Hurwitz condition applied to (\ref{J_comp}) is
\begin{equation}\label{1-RH}
{\textrm{tr}}J^C=b(1-3X^2)+c(1-3Y^2)<0.
\end{equation}
But for the assumptions (\ref{p13_assump}) it cannot be satisfied.
If only one of (\ref{p13_assump}) is satisfied, say the first one,
from the condition on the trace we obtain $b<-c(1-3Y^2)(1-3X^2)^{-1}$
and substituting into the determinant,
we have the estimate $\det J^C=b(1-3X^2)c(1-3Y^2)-a<-c^2(1-3Y^2)^2-a<0$ so that
the second Routh-Hurwitz condition is not satisfied. Hence the claim.

\vspace*{0.25cm}
\begin{remark}
Considering Figure \ref{scenario}, in the case of just one equilibrium, it must have at least one coordinate to the left (or below)
the one of the local maximum of the function. In the plot, it has the abscissa smaller than the one of the local maximum of the function
$Y_{[1]}(X)$. Thus when $E_3^C$ is unique, it must be unstable. For the case of three equilibria, evidently $E_4^C$ and $E_5^C$ have either
the abscissa ($E_4^C$) or the height ($E_5^C$)
satisfying the corresponding condition in (\ref{p13_assump}). Hence these two
equilibria must be unstable as well.
\end{remark}

\vspace*{0.25cm}

In case of three equilibria, we can show further the following result.

\vspace*{0.25cm}
\begin{proposition}
The equilibrium $E_3^C$ for which both the following conditions hold
\begin{equation}\label{p13_assump_opp}
X>\frac {\sqrt 3}3, \quad Y>\frac {\sqrt 3}3
\end{equation}
is stable.
\end{proposition}

\noindent {\textbf{Proof}}.
The Routh-Hurwitz condition (\ref{1-RH}) easily holds. The second one applied to (\ref{J_comp}) requires
$$
\det J^C=b(1-3X^2)c(1-3Y^2)-a>0.
$$
Observe that the slope of $Y_{[1]}(X)$ is negative at $X=1$. Hence for the abscissa of $E_3^C$ we must have $X_3<1$. Similarly $Y_3<1$,
using the slope of $X_{[2]}(Y)$ at $Y=1$. It follows that $b(1-3X^2)>-2b$, $c(1-3Y^2)>-2c$. Thus in turn $\det J^C>4bc-a$. Since we are in the case $a<bc$,
$\det J^C>0$ follows.

\vspace*{0.25cm}

\begin{remark}
There is thus a subcritical pitchfork bifurcation for which from the unstable $E_3^C$ three equilibria arise, with the equilibrium $E_3^C$ becoming
stable and the other ones being unstable.
\end{remark}

\vspace*{0.25cm}
\begin{remark}
No Hopf bifurcations arise in this model as they do not in the symbiotic one. Using the same technique as in the proof of
Proposition 14, the condition on the trace becomes an equality, so that by solving it
for $b$ we get $b=-c(1-3Y^2)(1-3X^2)^{-1}$. Substituting into the second Routh-Hurwitz condition $\det J^C>0$,
we obtain the contradiction $-c^2(1-3Y^2)^2-a>0$.
\end{remark}

\begin{figure}[!htb]
\begin{minipage}[b]{0.3\textwidth}
\includegraphics[width=1.1\textwidth]{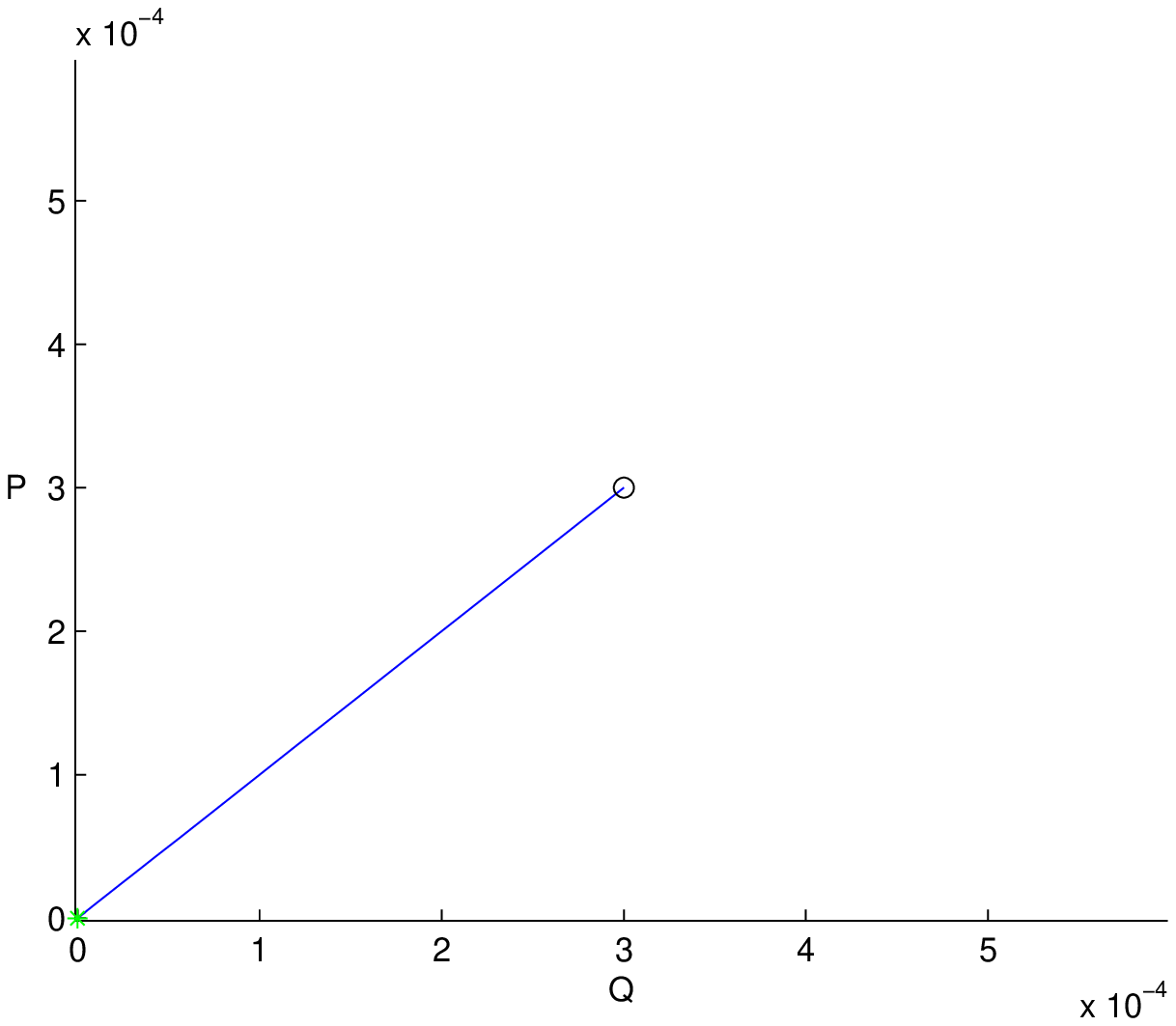}
\label{origin}
\end{minipage}
\begin{minipage}[b]{0.3\textwidth}
\includegraphics[width=1.1\textwidth]{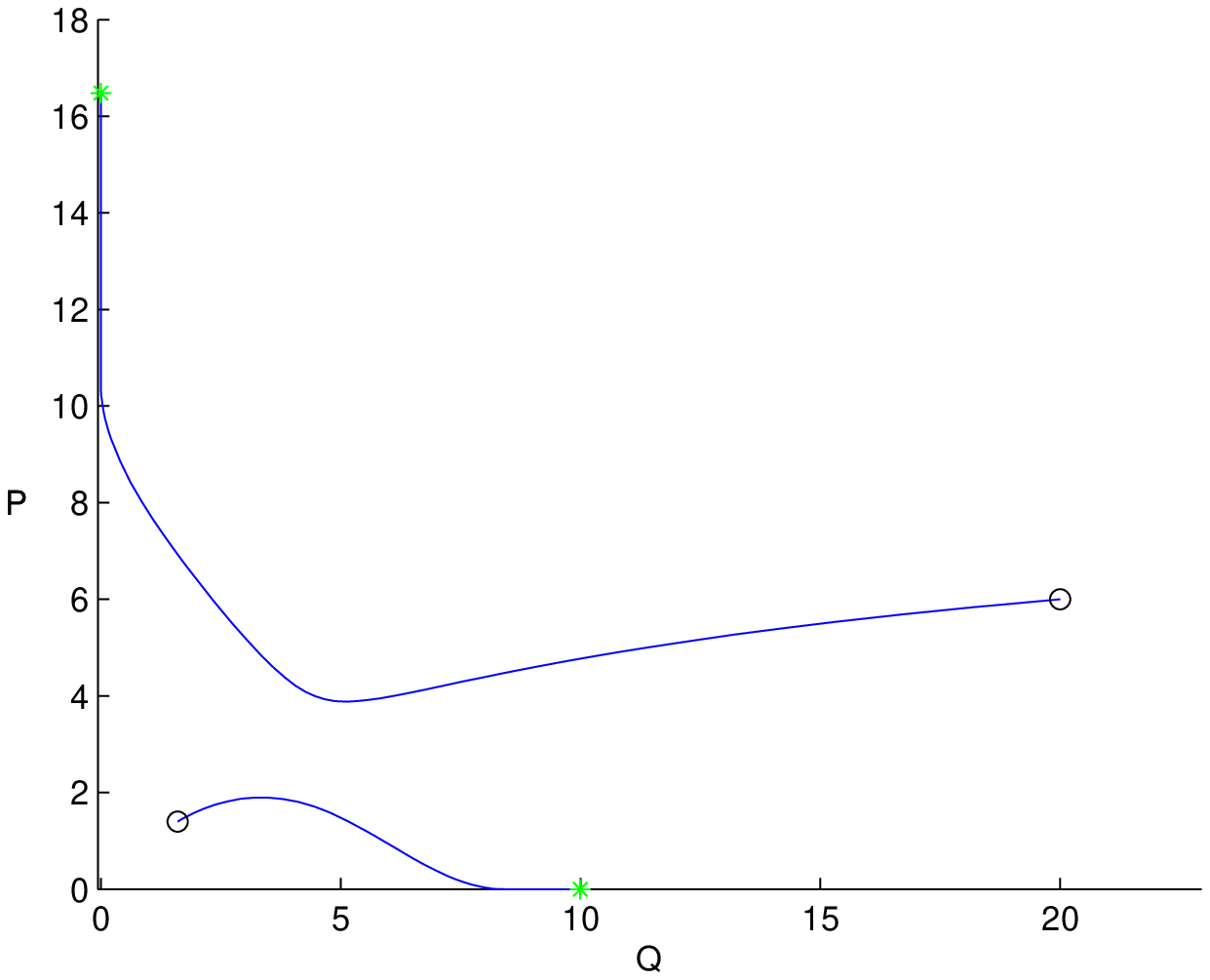}
\label{bistab}
\end{minipage}
\begin{minipage}[b]{0.3\textwidth}
\includegraphics[width=1.1\textwidth]{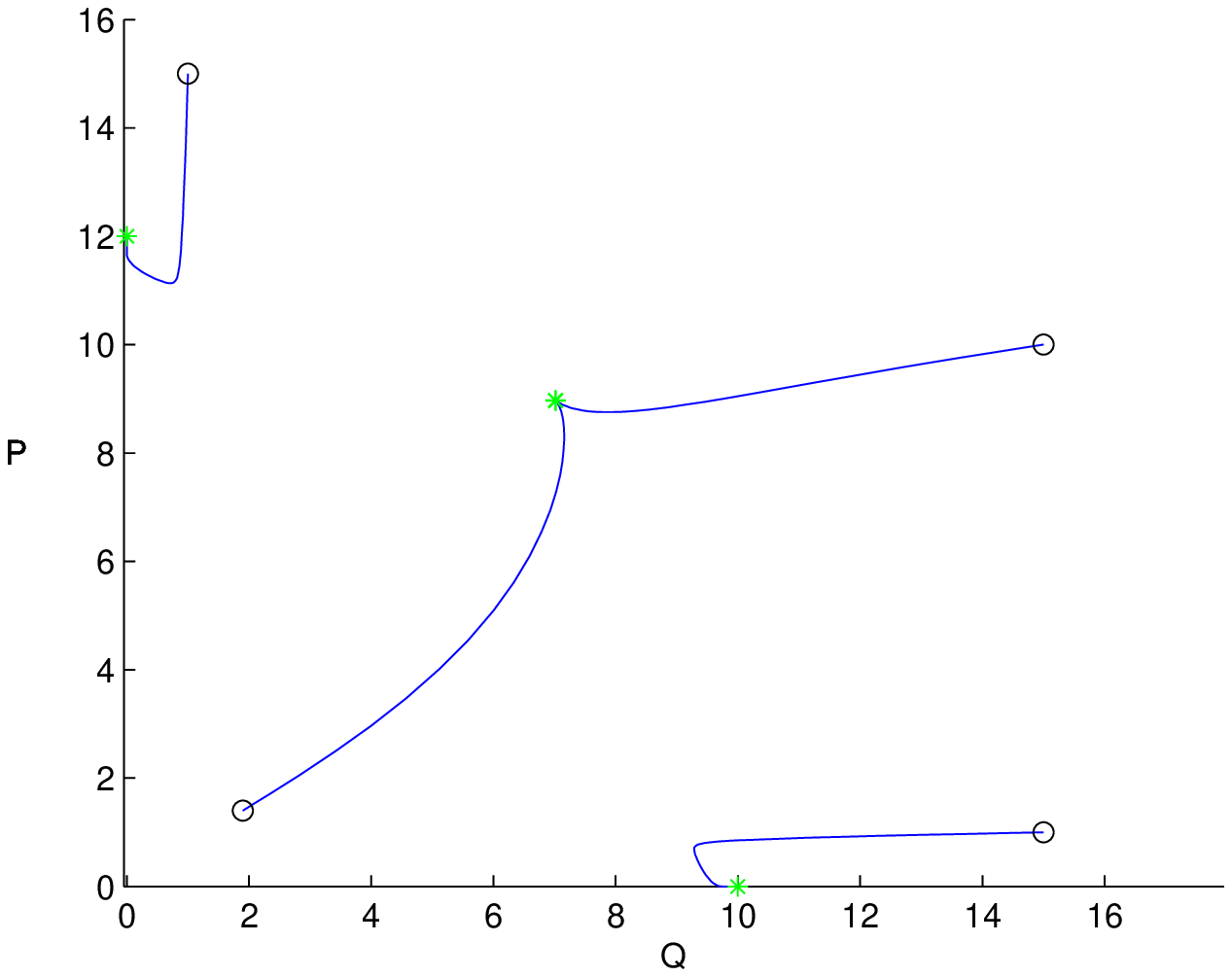}
\label{tristab}
\end{minipage}
\caption{The three possible populations behaviors. Left: the origin is stable, both populations get extinguished; achieved with parameter values
$a=0.75$, $b=0.0525$, $c=0.0525$, $r=2$, $m=2$, $p=33$, $q=33$, $K_p=4$, $K_q=3$.
Center: bistability and competitive exclusion, only one population survives; achieved with parameter values
$a=0.4052$, $b=1.1536$, $c=0.7814$, $r=0.8888$, $m=0.602$, $p=0.401$, $q=0.5998$, $K_p=16.5$, $K_q=10$.
Right: tristability, either one population only survives, or the other one, or both together; achieved with parameter values
$a=0.8993$, $b=3.4567$, $c=3.4523$, $r=0.7895$, $m=0.7885$, $p=0.225$, $q=0.2085$, $K_p=12$, $K_q=10$.
The green full dots represent the stable equilibria, the empty red circles are instead the initial conditions.}
\label{fig:tri}
\end{figure}

\vspace*{0.25cm}
In Figure \ref{fig:tri} we show the behavior of the two populations in the three possible cases.

\subsection{Comparison with classical competition model}

The classical competition model,
\begin{eqnarray}\label{comp_class}
\frac{dQ}{d\tau}=r\left( 1-\frac{Q}{K_Q}\right)Q-q PQ, \quad
\frac{dP}{d\tau}=m\left( 1-\frac{P}{K_P}\right)P-p PQ,
\end{eqnarray}
shows under suitable circumstances the competitive exclusion principle.
Thus, only one population survives, while the other one is wiped out. The system's outcome depends only
on its initial conditions, so that if the system has population values lying in the attracting set of
one of the equilibria, the dynamics will be drawn to it unless the environmental conditions, i.e. the parameters in the model, abruptely change.

Instead, we have found here that in presence of community behavior of both populations, the same occurs, but there is another possibility,
namely tristability. When the conditions arise, the coexistence equilibrium may be present together with the equilibria in which one population
vanishes. Therefore the system's outcome is once more determined by the initial conditions, but this time the phase plane is partitioned
into three basins of attractions, corresponding each to one of the possible equilibria. It would be interesting to compute explicitly
the boundary of each one of them. For this task an
extension of the algorithms presented in \cite{Rodi,CDRPV} would be needed.

\begin{figure}[!htb]
\begin{minipage}[b]{0.45\linewidth}
\centering
\includegraphics[width=\textwidth]{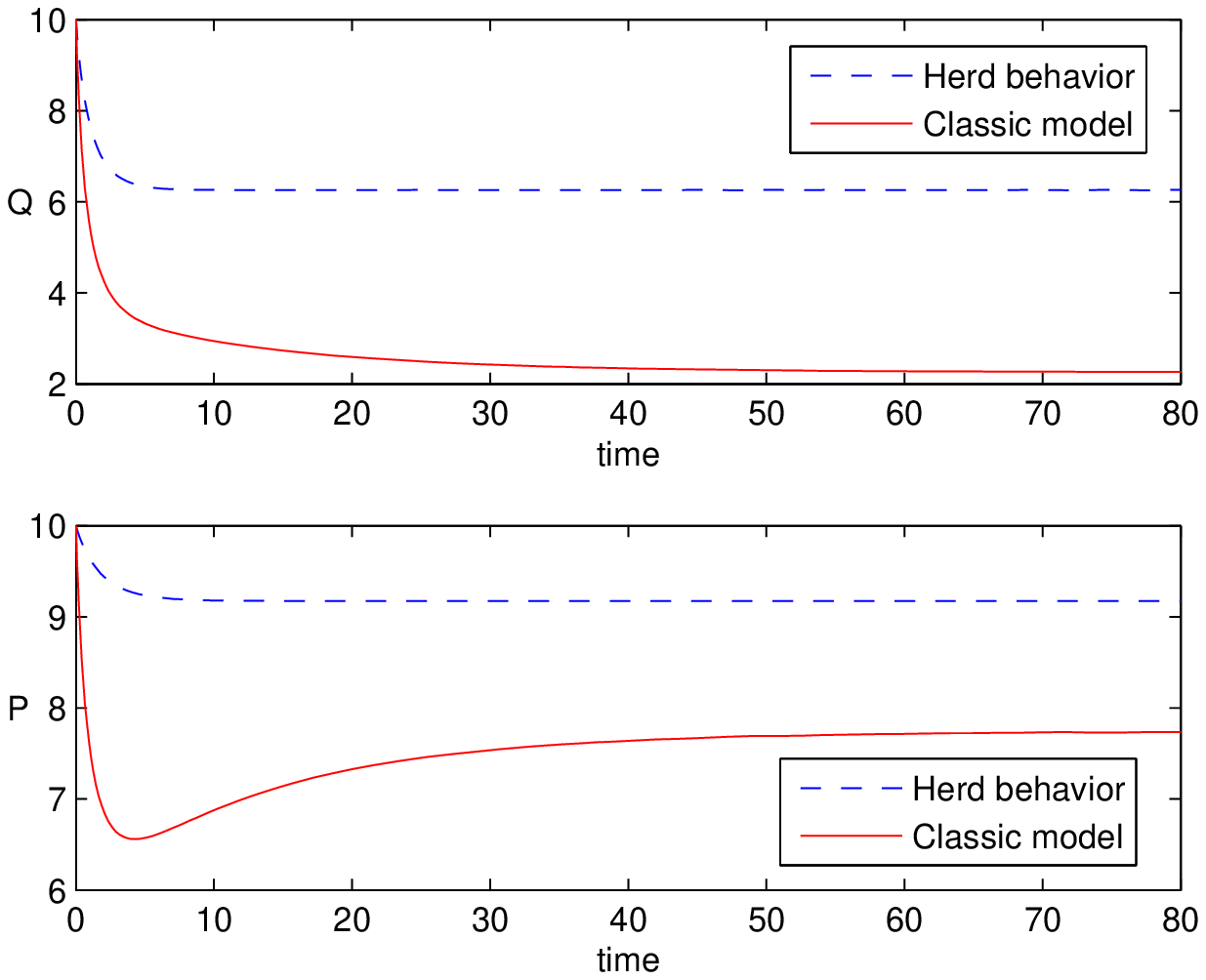}
\end{minipage}
\hspace{0.5cm}
\begin{minipage}[b]{0.45\linewidth}
\centering
\includegraphics[width=\textwidth]{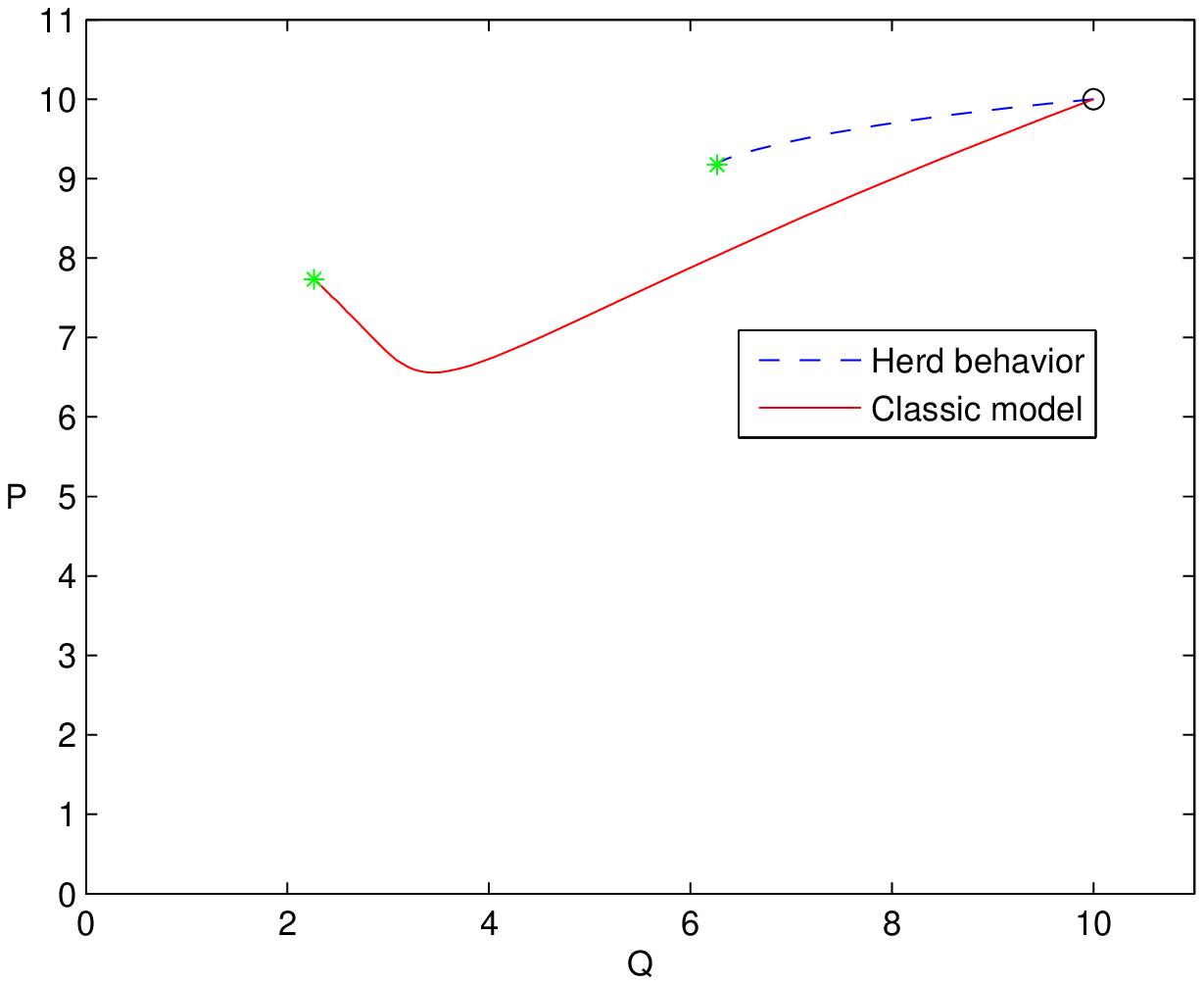}
\end{minipage}
\caption{Left: time series of the systems trajectories; Right: phase plane for classical and new symbiotic model with
$a=0.5$, $b=9.5618$, $c=5.9761$, $r=0.8$, $m=0.5$, $p=0.05$, $q=0.07$, $K_p=10$, $K_q=7$. The full green dots represent the equilibrium points.}
\label{compare_competition}
\end{figure}

We now compare the population levels when a coexistence equilibrium is stable in both classical and new model.
Considering the parameters $r=2$, $m=3$, $K_Q=6$, $K_P=8$, $q=0.2$ and $p=0.09$, with suitable initial conditions,
the behavior of the two models is shown in Figure \ref{compare_competition}. From the same initial conditions, trajectories of the two models
evolve toward different equilibria.
The population levels are thus higher in the herd model, at $Q^C=6.26$ and $P^C=9.17$ while for the classic model we find
$\widetilde Q^C=2.26$ and $\widetilde P^C=7.74$.
This is not surprising for the same reasons for which the opposite behavior occurs in the symbiotic models.
In herd models, only individuals at the outskirts meet individual of the other species. This means that individuals at the centre of
the flock here
receive less harm from the competition. On the contrary, in the classic model, individuals of the two populations are mixed together,
so that the whole populations are harmed by the competition.  

\section{Conclusions}

We have presented four models for non-classical population interactions, in that the populations involved in some way exhibit a socialized way of living.
This investigation completes the one undertaken in \cite{1}, in that all the situations that are possible in terms of individualistic or
gathering populations behavior are now analysed.
The models missing in \cite{1} are presented here:
we allow predators to hunt in packs, as well as
both intermingling populations to gather together, in the two cases of symbiosis and competition, so that they interact not on an individualistic basis,
but rather is some coordinate fashion.

The newly introduced symbiotic model on a qualitative basis behaves like the classical one. The populations settle always at the
coexistence equilibrium. Only, their levels are quantitatively smaller than in the classical case since the mutually beneficial
interactions in the new model are somewhat reduced.

For predator-prey interactions in presence of predators' pack hunting, we may have the prey behave in herds or individualistically.
The most prominent discrepancy between these two cases is the fact that both populations may disappear, under specific unfortunate conditions,
when the prey use a defensive coordinate strategy. This does not happen instead if they move loose in the environment, i.e. exhibit individualistic behavior,
since they attain a coexistence equilibrium.
This finding is quite
counterintuitive, because it could imply that the defensive mechanism is ineffective. But an interpretation could be provided,
since herds are more easily spotted by predators than individuals who can more easily hide in the terrain configuration.
Once the prey herds are completely wiped out, the predators also will disappear, since they are assumed not to be generalist, i.e.
their only food source is the prey under consideration.
Ecosystem extinction has also been rarely observed in the model without pack predation, \cite{PV}.
The system with prey herd behavior also shows limit cycles, i.e. the populations can coexist also through persistent oscillations,
not only at a stable equilibrium, which instead is the only possible system's outcome for the model with individualistic prey.
A similar result had been discovered earlier in case of individualistic predators hunting, \cite{1},
constituting the major difference between the prey group defense model with uncoordinated predation and the classical predator-prey system.
Finally, on the quantitative side, the coexistence population values for these two models with pack hunting differ, but without specific
informations on the parameter values it is not possible to assess which system will provide higher population values.

The competition system presented here allows again the extinction of both populations, under unfavorable circumstances, while this
never happens for the classical model.
Ecosystem disappearance occurs when (\ref{comp_disapp}) holds,
a condition that in the nondimensional model is equivalent to $a>bc$, as stated in Proposition 12.
When the competition system thrives, it does at higher levels for both populations than those achieved in the classical model.
Thus in this case populations coordinated behavior boosts their respective sizes, in case the system
parameters are in the range for which coexistence occurs.

But the major finding in this context of social behavior among all possible populations behavior is found for
the competition case. Indeed the system in suitable conditions
can show the phenomenon of competitive exclusion as the classical model does, but in addition we have discovered
that both populations can thrive, together with the situations predicted by the competitive exclusion principle. In other words, we
have found that the rather simple model (\ref{comp}) or (\ref{adimc}) may exhibit tristability, see once more the right picture in Figure \ref{fig:tri}.
This appears to be a novel and quite interesting finding further characterizing the systems with socialized behaviors.
The authors do not know of any other simple related model with such behavior.

%\vspace*{0.5cm}

%\section*{Acknowledgements}
%This paper is dedicated to the 70th birthday of M. Mackey.


\begin{thebibliography}{99}

\bibitem{1} V.~Ajraldi, M.~Pittavino, E.~Venturino, \textit{Modeling herd behavior in population systems},
Nonlinear Analysis: Real World Applications,  12 (2011) 2319-2338.

\bibitem{4} P.A.~Braza, \textit{Predator prey dynamics with square root functional responses},
Nonlinear Analysis: Real World Applications, 13 (2012) 1837-1843.~

\bibitem{B-VdD}  S. Busenberg, P. van den Driessche, (1990) Analysis of a
disease transmission model in a population with varying size, J. of
Math. Biology, 28, p. 257-270.

\bibitem{Samrat_NOVA} E.~Caccherano, S.~Chatterjee, L.~Costa Giani, L.~Il Grande, T.~Romano, G.~Visconti,
E.~Venturino, \textit{Models of symbiotic associations in food chains}, in
Symbiosis: Evolution, Biology and Ecological Effects, A.F.~Camis\~ao and C.C.~Pedroso (Editors),
Nova Science Publishers, Hauppauge, NY, 189-234, 2012.

\bibitem{CV} E.~Cagliero, E.~Venturino, \textit{Ecoepidemics with group defense and infected prey protected by herd},
Proceedings of the 12th International Conference on Computational and 
Mathematical Methods in Science and Engineering, CMMSE 2012,  J.~Vigo-Aguiar, A.P.~Buslaev, A.~Cordero, M.~Demiralp, I.P.~Hamilton,
E.~Jeannot, V.V.~Kozlov, M.T.~Monteiro, J.J.~Moreno, J.C.~Reboredo, P.~Schwerdtfeger, N.~Stollenwerk, J.R.~Torregrosa, E.~Venturino, J.~Whiteman (Editors)
1 (2012) 247-266.

\bibitem{Rodi} R.~Cavoretto, S.~Chaudhuri, A.~De Rossi, E.~Menduni,
F.~Moretti, M.C.~Rodi, E.~Venturino,
\textit{Approximation of Dynamical System's Separatrix Curves},
Numerical Analysis and Applied Mathematics ICNAAM 2011, T.~Simos, G.~Psihoylos, Ch.~Tsitouras, Z.~Anastassi (Editors),
AIP Conf.~Proc. 1389, 1220-1223 (2011); doi: 10.1063/1.3637836.

\bibitem{CDRPV} R.~Cavoretto, A.~De Rossi, E.~Perracchione, E.~Venturino,
\textit{Reconstruction of separatrix curves and surfaces in squirrels competition models with niche},
Proceedings of the 2013 International Conference on Computational and Mathematical Methods in Science and Engineering,
I.P.~Hamilton, J.~Vigo-Aguiar, H.~Hadeli, P.~Alonso, M.T.~De Bustos, M.~Demiralp, J.A.~Ferreira, A.Q.M.~Khaliq, J.a.~L\'opez-Ramos,
P.~Oliveira, J.C.~Reboredo, M.~Van Daele, E.~Venturino, J.~Whiteman, B.~Wade (Editors)
Almeria, Spain, June 24th-27th, 2013, v. 3, p. 400-411.

\bibitem{cosner} C.~Cosner, D.L.~DeAngelis, J.S.~Ault, D.B.~Olson, \textit{Effects of Spatial Grouping on the Functional
Response of Predators}, Theoretical Population Biology, 56, (1999) 65-75.

\bibitem{FW} H.I.~Freedman, G.~Wolkowitz,
\textit{Predator-prey systems with group defence: the paradox of enrichment revisited},
Bull. Math. Biol., 48 (1986) 493-508.

\bibitem{G-H}  L.Q. Gao, H.W. Hethcote, (1992), Disease transmission models
with density-dependent demographics, J. of Math. Biology, 30, p.
717-731.

\bibitem{GG} S.A.H. Geritz, M. Gyllenberg, Group defence and the predator's
functional response, Journal of Mathematical Biology 66, (2013) 705-717.

\bibitem{Giardina} I. Giardina,
Collective behavior in animal groups: theoretical models and empirical studies,
HFSP J. 2008 August; 2(4): 205-219, doi:  10.2976/1.2961038.

\bibitem{NOVA} M.~Haque, E.~Venturino, \textit{Mathematical models of diseases spreading in symbiotic communities},
in J.D.~Harris, P.L.~Brown (Editors), Wildlife: Destruction, Conservation and Biodiversity, NOVA Science
Publishers, New York, 2009, 135-179.

\bibitem{H00} Hethcote, H.~W., The mathematics of infectious diseases, SIAM Review 42 (2000) 599-653.

\bibitem{MPV} H.~Malchow, S.~Petrovskii, E.~Venturino, \textit{Spatiotemporal patterns in Ecology and Epidemiology},
CRC, Boca Raton, 2008.

\bibitem{ML-H}  J. Mena-Lorca, H.W. Hethcote, (1992), Dynamic models of
infectious diseases as regulator of population sizes, J. Math. Biology,
30, p. 693-716.

\bibitem{Sumpter} D.J.T Sumpter,
The principles of collective animal behaviour,
Phil. Trans. R. Soc. B 29, v. 361, n. 1465 (2006) 5-22,
doi: 10.1098/rstb.2005.1733 

\bibitem{EV_JBS} E. Venturino, A minimal model for ecoepidemics with group defense,
J. of Biological Systems 19(4), 763-785, 2011.

\bibitem{PV} E.~Venturino, S.~Petrovskii, \textit{Spatiotemporal Behavior of a Prey-Predator System with a Group Defense for Prey},
Ecological Complexity, 14 (2013) 37-47.

\end{thebibliography}
\end{document}